\documentclass[leqno,12pt]{article}
\usepackage{amsmath,amssymb,rotating,a4wide,color}
\usepackage{wasysym}

\usepackage[all,cmtip]{xy}

\usepackage{verbatim}

\newdir^{ (}{{}*!/-3pt/\dir^{(}}
\newdir_{ (}{{}*!/-3pt/\dir_{(}}

\entrymodifiers={+!!<0pt,\fontdimen22\textfont2>}

\def\sgn{\mathop{\rm sgn}\nolimits}
\def\bigtimes{\mathop{\raise-2pt\hbox{\huge$\times$}}}
\def\Ux{{\underline{x}}}
\def\GammaUx{\Gamma_{\kern-1pt\Ux}}
\def\OLGF{$\Phi$}
\def\tildew{{\widetilde{w}}}

\newbox\circbulletbox
\setbox\circbulletbox=\hbox{\,{\large$\circledcirc$}\kern-8.7pt\raise .6pt\hbox{$\bullet$}\,}

\let\le\leqslant
\let\ge\geqslant

\def\longland{\ \land\ }
\def\mycirc{{\kern1pt\circ\kern2pt}}

\def\Desc{\mathop{\rm Desc}\nolimits}

\def\Aut{\mathop{\rm Aut}\nolimits}

\def\deg{\mathop{\rm deg}\nolimits}

\def\Ker{\mathop{\rm Ker}\nolimits}

\def\ord{\mathop{\rm ord}\nolimits}

\let\phi\varphi
\let\theta\vartheta
\let\epsilon\varepsilon
\let\setminus\smallsetminus
\let\ndiv\nmid
\let\emptyset\varnothing

\newtheorem{Thm}{Theorem}[section]
\newtheorem{Prop}[Thm]{Proposition}
\newtheorem{Lem}[Thm]{Lemma}

\newtheorem{Sublem}[Thm]{Sublemma}
\newtheorem{Cor}[Thm]{Corollary}

\newtheorem{Def}[Thm]{Definition}

\newtheorem{Cond}[Thm]{Condition}
\newtheorem{Conds}[Thm]{Conditions}

\numberwithin{Thm}{section}
\def\UseTheoremCounterForNextEquation{\setcounter{equation}{\value{Thm}}\addtocounter{Thm}{1}}

\def\qed{{\hskip0pt\unskip\unskip\nobreak\hfil\penalty50
          \hskip1em\hbox{}\nobreak\hfil
           {$\square$}
          \parfillskip=0pt\finalhyphendemerits=0
          \par}\medskip}

\newenvironment{Proof}
               {\noindent{\bf Proof.}\ }
               {\qed}

               {\noindent{\bf Proof of #1.}\ }
               {\qed}

\newcommand{\BC}{{\mathbb{C}}}

\newcommand{\BF}{{\mathbb{F}}}

\newcommand{\BP}{{\mathbb{P}}}
\newcommand{\BQ}{{\mathbb{Q}}}

\newcommand{\BZ}{{\mathbb{Z}}}

\newbox\mybox
\def\arrover#1{\mathrel{
       \setbox\mybox=\hbox spread 1.4em
              {\hfil$\scriptstyle#1$\hfil}
       \vbox{\offinterlineskip\copy\mybox
             \hbox to\wd\mybox{\rightarrowfill}}}}

\def\larrover#1{\mathrel{
       \setbox\mybox=\hbox spread 1.4em
              {\hfil$\scriptstyle#1\vphantom{g}$\hfil}
       \vbox{\offinterlineskip\copy\mybox
             \hbox to\wd\mybox{\leftarrowfill}}}}

\def\ontoover#1{\mathrel{
       \setbox\mybox=\hbox spread 1.4em
              {\hfil$\scriptstyle#1\vphantom{g}$\hfil}
       \vbox{\offinterlineskip\copy\mybox
             \hbox to\wd\mybox{\rightarrowfill\hskip-2.8mm
                               $\rightarrow$}}}}
\def\leftontoover#1{\mathrel{
       \setbox\mybox=\hbox spread 1.4em
              {\hfil$\scriptstyle#1\vphantom{g}$\hfil}
       \vbox{\offinterlineskip\copy\mybox
             \hbox to\wd\mybox{$\leftarrow$\hskip-2.8mm
                               \leftarrowfill}}}}

\begin{document}

\title{Orbit length generating functions of\\
automorphisms of a rooted regular binary tree}

\author{Richard Pink\\[12pt]
\small Department of Mathematics \\[-3pt]
\small ETH Z\"urich\\[-3pt]
\small 8092 Z\"urich\\[-3pt]
\small Switzerland \\[-3pt]
\small pink@math.ethz.ch\\[12pt]}


\date{March 27, 2014}

\maketitle

\begin{abstract}
To every automorphism $w$ of an infinite rooted regular binary tree we associate a two variable generating function $\Phi_w$ that encodes information on the orbit structure of~$w$. We prove that this is a rational function if $w$ can be described by finitely many recursion relations of a particular form.
We show that this condition is satisfied for all elements of the discrete iterated monodromy group $\Gamma$ associated to a postcritically finite quadratic polynomial over~$\BC$. For such $\Gamma$ we also prove that there are only finitely many possibilities for the denominator of~$\Phi_w$, and we describe a procedure to determine their lowest common denominator.
\end{abstract}

{\renewcommand{\thefootnote}{}
}

\newpage
\tableofcontents
\newpage


\section{Introduction}
\label{Intro}

Let $T$ be an infinite rooted regular binary tree. To any automorphism $w$ of $T$ we associate the power series
$$\Phi_w \ = \sum_{n\ge m\ge0} o_{m,n}(w)X^mY^n \ \in\ 1+Y\BZ[[X,Y]]$$
where $o_{m,n}(w)$ is the number of orbits of $w$ of length $2^m$ on the set of vertices of level $n$ of~$T$. This \emph{orbit length generating function} encodes some information, but in general not all, about the conjugacy class of~$w$. 

\medskip
The use of this construction stems from its behavior with respect to recursion relations. Specifically, assume that we are given an isomorphism from $T$ to each half subtree of $T$ obtained after deleting the root. Then for any $u,v\in\Aut(T)$ there is a unique element $(u,v)\in\Aut(T)$ which acts on the two half subtrees by $u$ and~$v$, respectively. Also, let $\sigma\in\Aut(T)$ denote the involution which interchanges the two half subtrees. Then one easily shows that
\UseTheoremCounterForNextEquation
\begin{equation}\label{PhiRecIntro}
\begin{array}{rl}
\Phi_{(u,v)} &\!\!=\ 1+Y\Phi_u+Y\Phi_v \quad\hbox{and}\\
\Phi_{(u,v)\,\sigma} &\!\!=\ 1+XY\Phi_{uv}.
\end{array}
\end{equation}
These relations are particularly useful for automorphisms that are themselves described by recursion relations. Many such automorphisms can be described abstractly, but they arise most notably as elements of iterated monodromy groups of quadratic morphisms.

\medskip
To apply the recursion relations we say that an element $w\in\Aut(T)$ is \emph{\OLGF-finite} if repeated application of the rules $(u,v)\rightsquigarrow u, v$ and $(u,v)\,\sigma \rightsquigarrow uv$ beginning with $w$ leads to only finitely many elements of $\Aut(T)$. We say that $w$ is \emph{\OLGF-irreducible} if, in addition, the rules eventually lead back to~$w$. 
Using (\ref{PhiRecIntro}) it is not hard to prove that for any \OLGF-finite element $w$ the power series $\Phi_w$ is the expansion of a rational function in $X$ and~$Y$, and that new factors in the denominator arise only for \OLGF-irreducible elements (Theorems \ref{OLGFRat} and \ref{OLGFLinComb}).

\medskip
Roughly speaking an element $w$ is \OLGF-finite if and only if it results from finitely many recursion relations of some particular form. This condition is probably quite restrictive. In fact, for the automorphism defined by the relatively easy looking recursion relation $b=(b,b\sigma)\,\sigma$ we explicitly calculate $\Phi_b$ and show that it is not a rational function (see Section \ref{NonRat}).

\medskip
By contrast, fix two integers $r>s\ge0$ and consider a tuple $\Ux=(x_2,\ldots,x_r)$ with entries in $\{0,1\}$. To this data Bartholdi and Nekrashevych \cite{Bartholdi-Nekrashevych-2008} have associated a certain subgroup $\GammaUx \subset\Aut(T)$ by explicit recursion relations for $r$ generators. They have shown that the iterated monodromy group of any quadratic polynomial in one variable over~$\BC$ with a finite postcritical orbit of size $r$ and eventual period $r-s$ is conjugate to $\GammaUx$ for some~$\Ux$.

\medskip
{}From the recursion relations of the generators alone we deduce with modest effort that all elements of $\GammaUx$ are \OLGF-finite (Propositions \ref{0OLGFFin} and \ref{1OLGFFin}). 
More surprisingly, and with much more work, we prove that the rational functions $\Phi_w$ for all $w\in\GammaUx$ possess a common denominator that depends only on~$\GammaUx$. We describe a common denominator explicitly and characterize the unique lowest common denominator $D_\Ux$ in terms of a finite combinatorial problem concerning the data $r$, $s$, and~$\Ux$ (Theorems \ref{0MainThm} and \ref{1MainThm}).

\medskip
These results rely on a detailed analysis of the \OLGF-irreducible elements in~$\GammaUx$. While there are infinitely many of them, we show that they lie in an explicit finite collection of conjugacy classes of $\Aut(T)$ (Propositions \ref{0IrredConj} and \ref{1IrredConj}). As $\Phi_w$ is invariant under conjugacy, this implies the existence of some common denominator of $\Phi_w$ for all $w\in\GammaUx$. The characterization of the lowest common denominator requires additional effort.
 
\medskip
The results of this article lead to a number of interesting questions and open problems. Among these are:
\begin{enumerate}
\item[$\bullet$] 
By Nekrashevych \cite[Thm.\;6.4.4]{Nekrashevych-2005} the iterated monodromy group of any postcritically finite rational function over $\BC$ is \emph{contracting} in the sense of \cite[Def.\;2.11.1]{Nekrashevych-2005}. Our notion of \OLGF-finiteness is similar, but not equivalent. Are there deeper connections? Also, is there a relation with the notions of \emph{finite state}, \emph{bounded}, and/or \emph{finitary} automorphisms from Bartholdi-Nekrashevych \cite[\S2.4]{Bartholdi-Nekrashevych-2008}?

\item[$\bullet$] 
When $\Gamma_\Ux$ is the iterated monodromy group of a postcritically finite quadratic polynomial over~$\BC$, what do \OLGF-finiteness and the rationality of $\Phi_w$ mean geometrically? What is the geometric meaning of the numerator and denominator of~$\Phi_w$, and of the lowest common denominator of all~$\Phi_w$?

\item[$\bullet$] 
Prove \OLGF-finiteness and rationality and describe the denominators directly for iterated monodromy groups, without using their classification \`a la Bartholdi-Nekrashevych~\cite{Bartholdi-Nekrashevych-2008}, perhaps in a way similar to how the group theoretic contracting property is deduced from geometric facts.

\item[$\bullet$] 
Based on the polynomial case, we conjecture that for the iterated monodromy group $\Gamma$ of any postcritically finite quadratic morphism $\BP^1_\BC\to\BP^1_\BC$, the $\Phi_w$ for all $w\in\Gamma$ are rational and possess a common denominator.

\item[$\bullet$] 
Define orbit length generating functions for automorphisms of an infinite $d$-regular rooted tree for an arbitrary, possibly composite, integer $d\ge2$ and extend the present results accordingly.

\item[$\bullet$] 
As part of our analysis we prove that any \OLGF-irreducible element of $\GammaUx$ is conjugate under $\Aut(T)$ to some \OLGF-irreducible element of $\Gamma_{(0,\ldots,0)}$ with the same pair $(r,s)$ (Propositions \ref{0IrredConjToSplit} and \ref{1IrredConjToSplit}). Is the analogue true for non-\OLGF-irreducible elements? 

\item[$\bullet$] 
The role of $\Gamma_{(0,\ldots,0)}$ as a receptacle for conjugacy classes resembles the way that a quasi-split connected reductive group $G$ over a field $K$ possesses $K$-rational elements corresponding to the conjugacy classes of all $K$-rational elements of all inner forms of~$G$. Is there a similar sense in which $\Gamma_{(0,\ldots,0)}$ is a `quasi-split inner form of~$\GammaUx$'?

\item[$\bullet$] 
Our results show that the lowest common denominator $D_\Ux$ varies with~$\Ux$ and that it is largest when $\Ux=(0,\ldots,0)$ or $(1,\ldots,1)$. Since $D_\Ux$ depends only on the $\Aut(T)$-conjugacy class of~$\GammaUx$, it can help distinguish some of these conjugacy classes, especially from the conjugacy class of $\Gamma_{(0,\ldots,0)}$.
However, there are still many different tuples $\Ux$ with the same~$D_\Ux$. Can these groups be distinguished using the precise form of~$\Phi_w$, or using the conjugacy classes of non-\OLGF-irreducible elements?

\item[$\bullet$] 
Do the orbit length generating functions also distinguish the Grigorchuk group from its `twisted twin' of Bartholdi-Siegenthaler \cite{Bartholdi-Siegenthaler-2010}?

\item[$\bullet$] 
If the iterated monodromy groups associated to two postcritically finite quadratic polynomials over $\BC$ are conjugate in $\Aut(T)$, does it follow that the polynomials are equivalent under an affine linear transformation and/or complex conjugation?

\item[$\bullet$] 
Determine all subsets $J\subset\{1,\ldots,r\}$ satisfying Condition \ref{0JCond}, 
respectively Conditions \ref{1JConds}. 
Give a direct formula for the lowest common denominator $D_\Ux$ instead of a finite algorithm.

\item[$\bullet$] 
Our original motivation was to understand the action of Frobenius elements associated to quadratic morphisms defined over finite or finitely generated fields of characteristic $\not=2$. However, preliminary sample calculations suggest that their orbit length generating functions behave differently from those of the discrete groups studied in the present article. Nevertheless this question should be studied further, maybe in connection with the approach of Boston-Jones \cite{JonesBoston2012}.
\end{enumerate}

 \newpage


\section{General definitions and results}
\label{Gen}


\subsection{Notation}
\label{D1}

Let $T$ be the infinite tree whose vertices are the finite words over the alphabet $\{0,1\}$ and where each vertex $t$ is connnected by an edge to the vertices $t0$ and $t1$. The empty word is called the \emph{root of~$T$}, making $T$ an infinite rooted regular binary tree. 

\medskip
Let $W$ denote the automorphism group of~$T$. For any elements $u,v\in W$ we let $(u,v)$ denote the element of $W$ defined by $t0\mapsto u(t)0$ and $t1\mapsto v(t)1$ for any word~$t$. 
This defines an isomorphism from $W\times W$ to the subgroup of $W$ that fixes the vertices $0$ and~$1$. We identify $W\times W$ with its image.
Let $\sigma\in W$ denote the element of order $2$ defined by $t0\mapsto t1\mapsto t1$ for any word~$t$, and let $\langle\sigma\rangle$ be the subgroup of $W$ generated by it. Then $W$ is the semidirect product $W = (W\times W)\rtimes\langle\sigma\rangle$.


\medskip
For any integer $n\ge0$, the \emph{level $n$} of~$T$ is the set of vertices at distance $n$ from the root, i.e., the set of $2^n$ words of length~$n$. Any element $w\in W$ fixes the root and thus permutes the level~$n$. We let $\sgn_n(w)$ denote the sign of the induced permutation of the level~$n$. 
Then $\sgn_1(\sigma)=-1$ and $\sgn_n(\sigma)=1$ for all $n\not=1$, and for any $u,v\in W$ we have $\sgn_{n+1}((u,v)) = \sgn_n(u)\cdot\sgn_n(v)$.

\medskip
For any $n\ge0$ let $T_n$ denote the finite subtree obtained by cutting off $T$ at level~$n$. The automorphism group of $T_n$ is a certain iterated wreath product of the group of two elements with itself and therefore a finite $2$-group. Thus for any $w\in W$, any orbit of $w$ on level $n$ has length $2^m$ for some integer $0\le m\le n$. The root of $T$ is the unique vertex on level $0$ and constitutes an orbit of length~$1$.


\subsection{Orbit length generating functions}
\label{D2}

\begin{Def}\label{OLGFDef}
The \emph{orbit length generating function of} $w\in W$ is the power series 
$$\Phi_w \ = \sum_{n\ge m\ge0} o_{m,n}(w)X^mY^n \ \in\ 1+Y\BZ[[X,Y]]$$
where $o_{m,n}(w)$ is the number of orbits of $w$ of length $2^m$ on level~$n$.
\end{Def}

\begin{Lem}\label{OLGFProps}
For any element $w\in W$ we have:
\begin{enumerate}
\item[(a)] $\Phi_w$ depends only on the $W$-conjugacy class of~$w$.
\item[(b)] $\Phi_{w^k}=\Phi_w$ for any odd integer $k$.
\item[(c)] $\displaystyle \Phi_{w^2}(X,Y) = \Phi_w(0,Y) + 2\cdot\frac{\Phi_w(X,Y)-\Phi_w(0,Y)}{X}$.
\end{enumerate}
\end{Lem}

\begin{Proof}
Assertions (a) and (b) follow from the fact that the orbit lengths remain the same. Next, any fixed point of $w$ remains a fixed point of~$w^2$, and any orbit of length $2^{m+1}>1$ of $w$ splits into two orbits of length $2^m$ of~$w^2$. Thus $o_{0,n}(w^2) = o_{0,n}(w)+2o_{1,n}(w)$, and $o_{m,n}(w^2) = 2o_{m+1,n}(w)$ whenever $m>0$. This implies (c).
\end{Proof}

\begin{Prop}\label{OLGFRecRels}
For any elements $u,v\in W$ we have
\begin{eqnarray*}
\Phi_{(u,v)} &\!\!=\!\!& 1+Y\Phi_u+Y\Phi_v, \\
\Phi_{(u,v)\,\sigma} &\!\!=\!\!& 1+XY\Phi_{uv}.
\end{eqnarray*}
\end{Prop}

\begin{Proof}
By the definition of $(u,v)$, its orbits on level $n+1$ are obtained from the orbits of $u$ on level $n$ by appending the letter $0$ to each word
and from the orbits of $v$ on level $n$ by appending the letter $1$ each word. Thus $o_{m,n+1}((u,v)) = o_{m,n}(u)+o_{m,n}(v)$, which implies the first formula.

The other element $(u,v)\,\sigma$ fixes the root, but changes the last letter of every word of length $>0$. Thus its orbits of length $2^{m+1}$ are in bijection with the orbits of length $2^m$ of $(u,v)\,\sigma\,(u,v)\,\sigma = (uv,vu)$ on the set of words ending in~$0$. By the definition of $(uv,vu)$ the latter are obtained from the orbits of $uv$ of length $2^m$ by appending the letter $0$ to each word. Thus $o_{m+1,n+1}((u,v)\,\sigma) = o_{m,n}(uv)$, which implies the second formula.
\end{Proof}

The recursion relations in Proposition \ref{OLGFRecRels} are the main tools for calculating~$\Phi_w$. To formalize their use we introduce the following ad hoc terminology.


\subsection{Finiteness}
\label{D3}

\begin{Def}\label{DescDef}
The \emph{first descendants of} an element $w\in W$ are the elements $u$ and $v$ if $w=(u,v)$, respectively $uv$ alone if $w=(u,v)\,\sigma$. For any $n\ge1$, the first descendants of all $n^{\rm th}$ descendants of $w$ are the \emph{$(n+1)^{\rm st}$ descendants of~$w$}. The $n^{\rm th}$ descendants of $w$ for all $n\ge1$ are the \emph{descendants of~$w$}. The 
set of all descendants of $w$ is denoted $\Desc(w)$. 
\end{Def}

Thus $\Desc(w)$ is the set of elements of $W$ encountered on repeatedly applying the recursion relations \ref{OLGFRecRels}.

\begin{Def}\label{IrredDef}
\begin{enumerate}
\item[(a)] An element $w\in W$ is called \emph{\OLGF-finite} if $\Desc(w)$ is finite.
\item[(b)] An element $w\in W$ is called \emph{\OLGF-irreducible} if $\Desc(w)$ is finite and $w\in\Desc(w)$.
\end{enumerate}
\end{Def}


As a direct consequence of the definition we have:

\begin{Prop}\label{FinDesc}
For any $w'\in\Desc(w)$ we have $\Desc(w')\subset\Desc(w)$.
In particular, any descendant of a \OLGF-finite element is \OLGF-finite.
\end{Prop}


\subsection{Rationality}
\label{D4}

\begin{Thm}\label{OLGFRat}
If $w\in W$ is \OLGF-finite, then $\Phi_w$ is the power series expansion of a rational function in $X$ and $Y$ with denominator in $1+Y\BZ[X,Y]$.
\end{Thm}

\begin{Proof}
Write $\{w\}\cup\Desc(w) = \{w_1,\ldots,w_r\}$. Then Propositions \ref{OLGFRecRels} and \ref{FinDesc} imply that for any $1\le i\le r$ there exist $1\le j,k\le r$ such that $\Phi_{w_i} = 1+Y\Phi_{w_j}+Y\Phi_{w_k}$ or $\Phi_{w_i} = 1+XY\Phi_{w_j}$. 
In particular we can write $\Phi_{w_i} = 1+ \sum_{j=1}^r Ya_{i,j}\Phi_{w_j}$ for certain $a_{i,j}\in\BZ[X]$. In terms of the column vectors $f := (\Phi_{w_i})_{i=1}^r$ and $e := (1)_{i=1}^r$ and the matrix $A := (a_{i,j})_{i,j=1}^r$ this means that $f=e+YAf$. This in turn is equivalent to $(I-YA)f=e$, where $I$ denotes the identity matrix. The determinant $D := \det(I-YA)$ lies in $1+Y\BZ[X,Y]$ and is therefore invertible in $\BZ[[X,Y]]$, and the coefficients of $(I-YA)^{-1}$ lie in $D^{-1}\BZ[X,Y]$. Thus the coefficients of $f=(I-YA)^{-1}e$ lie in $D^{-1}\BZ[X,Y]$, and hence so does~$\Phi_w$, as desired.
\end{Proof}

\begin{Thm}\label{OLGFLinComb}
If $w\in W$ is \OLGF-finite, then $\Phi_w$ is a $\BZ[X,Y]$-linear combination of the $\Phi_{w'}$ for all \OLGF-irreducible $w'\in\Desc(w)$.
\end{Thm}

\begin{Proof}
By induction on the cardinality of $\{w\}\cup\Desc(w)$ we may assume that the assertion holds for all \OLGF-finite elements $w'\in W$ with $|\{w'\}\cup\Desc(w')| < |\{w\}\cup\Desc(w)|$. If $w$ is \OLGF-irreducible, there is nothing to prove. So assume that $w$ is not \OLGF-irreducible. 
Then for any $w'\in\Desc(w)$ we have $w\not\in\Desc(w')\subset\Desc(w)$. Thus $\{w'\}\cup\Desc(w')$ is a proper subset of $\{w\}\cup\Desc(w)$, and so by the induction hypothesis the assertion already holds for~$w'$. In particular, in the case $w=(u,v)$ the assertion holds for $u$ and~$v$, and in the case $w=(u,v)\,\sigma$ the assertion holds for $uv$. Thus with the recursion relations from Proposition \ref{OLGFRecRels} the assertion follows for~$w$, as desired.
\end{Proof}


\subsection{Examples}
\label{D5}

Now we do some simple examples. First, the identity element $1\in W$ is equal to $(1,1)$ and therefore \OLGF-irreducible. With Proposition \ref{OLGFRecRels} we find that $\Phi_1 = 1+2Y\Phi_1$ and so 
\UseTheoremCounterForNextEquation
\begin{equation}\label{F1}
\Phi_1\ =\ \frac{1}{1-2Y}.
\end{equation}
Next $\sigma=(1,1)\,\sigma$ has the unique descendant~$1$. Thus it is \OLGF-finite but not \OLGF-irreducible, and from (\ref{F1}) and Proposition \ref{OLGFRecRels} we deduce that
\UseTheoremCounterForNextEquation
\begin{equation}\label{Fsigma}
\Phi_\sigma \ =\ 1+\frac{XY}{1-2Y}.
\end{equation}
Next the \emph{standard odometer} is the element $a\in W$ defined by the recursion relation $a = (a,1)\,\sigma$. Thus it is \OLGF-irreducible, and from Proposition \ref{OLGFRecRels} we deduce that $\Phi_a = 1+XY\Phi_a$ and hence 
\UseTheoremCounterForNextEquation
\begin{equation}\label{FOdo}
\Phi_a \ =\ \frac{1}{1-XY}.
\end{equation}
Also, for any odd integer $k=2\ell+1$ the element $a^k = (a^{\ell+1},a^\ell)\sigma$ is again \OLGF-irreducible and has $\Phi_{a^k}=\Phi_a$ by Lemma \ref{OLGFProps} (b). In fact, one easily shows that any odd power of any \OLGF-irreducible element is \OLGF-irreducible.

\medskip
On the other hand, not all elements of $W$ that are described by finitely many recursion relations have rational orbit length generating functions, as the example in the next section shows. 
Also, rationality is rare in the following sense. Recall that as a profinite group $W$ has a unique Haar measure with total volume~$1$.

\begin{Prop}\label{MeasureZero}
The set of elements $w\in W$ with $\Phi_w$ rational has measure zero.
\end{Prop}

\begin{Proof}
As there are only countably many rational functions with coefficients in~$\BZ$, it suffices to prove that for any fixed $\Phi\in1+Y\BZ[[X,Y]]$, the set $S$ of all $w\in W$ with $\Phi_w=\Phi$ has measure zero. But $\Phi_w$ determines $\sgn_n(w)$ for all $n\ge1$, and so $S$ is contained in a single coset of the subgroup $\bigcap_{n\ge1}\Ker(\sgn_n)$ of~$W$. This is a closed subgroup of infinite index and therefore of measure zero; hence $S$ has measure zero, as desired.
\end{Proof}


\subsection{Variant}
\label{D6}

Some calculations become easier with the following slight variant of~$\Phi_w$ obtained by `removing trivial poles and zeros':

\begin{Prop}\label{OLGFRecRels21}
For any $w\in W$ there exists a unique $\Psi_w\in Y\BZ[[X,Y]]$ with 
$$\Phi_w \ =\ \frac{1}{1-2Y} \;+\; \frac{X-2}{1-2Y} \cdot \Psi_w.$$
\end{Prop}

\begin{Proof}
The term $o_{m,n}(w)$ in Definition \ref{OLGFDef} is the number of orbits of $w$ of length $2^m$ on level~$n$. Since the total number of vertices on level $n$ is~$2^n$, this implies that
$$\Phi_w(2,Y) \ =\ \sum_{n\ge m\ge0} o_{m,n}(w)2^mY^n
\ =\ \sum_{n\ge0}\; 2^nY^n
\ =\ \frac{1}{1-2Y}.$$
Thus $\Phi_w-\frac{1}{1-2Y}$ is divisible by~$X-2$, and the decomposition follows.
\end{Proof}

\begin{Prop}\label{OLGFRecRels22}
For any elements $u,v\in W$ we have
\begin{eqnarray*}
\Psi_{(u,v)} &\!\!=\!\!& Y\Psi_u+Y\Psi_v, \\
\Psi_{(u,v)\,\sigma} &\!\!=\!\!& Y+XY\Psi_{uv}.
\end{eqnarray*}
\end{Prop}

\begin{Proof}
Direct consequence of Proposition \ref{OLGFRecRels}.
\end{Proof}

\medskip
For example, the formulas (\ref{F1}) and (\ref{Fsigma}) and (\ref{FOdo}) correspond to:
\begin{eqnarray}
\UseTheoremCounterForNextEquation
\label{F12}
\Psi_1 &\!\!=\!\!& 0, \\[5pt]
\UseTheoremCounterForNextEquation
\label{Fsigma2}
\Psi_\sigma &\!\!=\!\!& Y, \\
\UseTheoremCounterForNextEquation
\label{FOdo2}
\Psi_a &\!\!=\!\!& \frac{Y}{1-XY}.
\end{eqnarray}

%
%
%
%
%
%
%
%

 \newpage


\section{A non-rational orbit length generating function}
\label{NonRat}

In this section we study the element $b\in W$ defined by the recursion relation
\UseTheoremCounterForNextEquation
\begin{equation}\label{bRec}
b \ =\ (b,b\sigma)\sigma.
\end{equation}
We will explicitly calculate $\Phi_b$ and show that it is not a rational function. This implies that the description of elements of $W$ by finitely many recursion relations does not guarantee that their orbit length generating functions are rational.


\subsection{Preparations}
\label{E1}

First note that the power $\sigma^p$ for $p\in\BZ$ depends only on $p\bmod2$ and can therefore be defined for any $p\in\BF_2$.
Thus to any integer $r\ge1$ and any polynomial $P(T) = \sum p_iT^i \in\BF_2[T]$ of degree $<2^r$ we can associate the element
\UseTheoremCounterForNextEquation
\begin{equation}\label{bLem0}
w_{r,P}\ :=\ b\,\sigma^{p_0}\,b\,\sigma^{p_1}\cdots b\,\sigma^{p_{2^r-1}} \ \in\ W.
\end{equation}
To any such $r$ and $P$ we also associate 
\begin{eqnarray*}
Q(T) &\!\!:=\!\!& \frac{P(T)\cdot T-P(1)\cdot T^{2^r}}{T-1} + T\cdot(T-1)^{2^r-2},
\\[3pt]
R(T) &\!\!:=\!\!& Q(T) + (T-1)^{2^r-1},
\qquad\hbox{and}\\[3pt]
S(T) &\!\!:=\!\!& Q(T) + T^{2^r}\cdot R(T),
\end{eqnarray*}
which are again polynomials in $\BF_2[T]$ of respective degrees $<2^r$, $<2^r$, and $<2^{r+1}$.

\begin{Lem}\label{bLem1}
In this situation $w_{r,P} = (w_{r,Q},w_{r,R})\,\sigma^{P(1)}$
and $w_{r,Q}w_{r,R}=w_{r+1,S}$.
\end{Lem}

\begin{Proof}
Set $q_i= \sum_{j=0}^{i-1}\,(p_j-1) \in \BF_2$ for all $0\le i\le 2^r$. 
Then $q_0=0$ and $p_i=1-q_i+q_{i+1}$ for all $0\le i<2^r$, and hence
$$w_{r,P}\ =\ (\sigma^{q_0}\,b\,\sigma^{1-q_0})\cdot
(\sigma^{q_1}\,b\,\sigma^{1-q_1})\cdots
(\sigma^{q_{2^r-1}}\,b\,\sigma^{1-q_{2^r-1}})\cdot\sigma^{q_{2^r}}.$$
Here $q_{2^r} = \sum_{j=0}^{2^r-1}\,(p_j-1) = \sum_{j=0}^{2^r-1}p_j - 2^r = P(1)$ because $r\ge1$. Also, for any $q\in\BF_2$ we have
$$\sigma^q\,b\,\sigma^{1-q}
\ =\ \sigma^q\,(b,b\sigma)\,\sigma^{-q}
\ =\ \biggl\{\!\begin{array}{ll}
(b,b\sigma) & \hbox{if $q=0$}\\[3pt]
(b\sigma,b) & \hbox{if $q=1$}
\end{array}\!\biggr\}
\ =\ (b\,\sigma^q,b\,\sigma^{q+1}).$$
Therefore
$$\begin{array}{rl}
w_{r,P} &=\ (b\,\sigma^{q_0},b\,\sigma^{q_0+1})\cdots
(b\,\sigma^{q_{2^r-1}},b\,\sigma^{q_{2^r-1}+1})
\cdot\sigma^{P(1)} \\[3pt]
&=\ \bigl( b\,\sigma^{q_0} \cdots b\,\sigma^{q_{2^r-1}},
b\,\sigma^{q_0+1}\cdots b\,\sigma^{q_{2^r-1}+1}\bigr)
\cdot\sigma^{P(1)}.
\end{array}$$
Thus with $Q(T) := \sum_{i=0}^{2^r-1} q_iT^i$ and $R(T) := \sum_{i=0}^{2^r-1} (q_i+1)T^i$ we deduce that $w_{r,P}=(w_{r,Q},w_{r,R})\,\sigma^{P(1)}$. A direct calculation which we leave to the reader shows that $Q(T)$ and $R(T)$ are given by the indicated formulas. Finally, the formula $w_{r,Q}w_{r,R}=w_{r+1,S}$ follows directly on expanding both sides.
\end{Proof}

As usual, for any polynomial $f\in\BF_2[T]$ we let $\ord_{T-1}(f)$ denote the supremum of the set of integers $d$ such that $(T-1)^d$ divides~$f$.

\begin{Lem}\label{bLem2}
If $0<\ord_{T-1}(P) < 2^r-1$, then $\ord_{T-1}(Q) = \ord_{T-1}(R) = {\ord_{T-1}(P)-1}$. Moreover, we always have $\ord_{T-1}(S) = 2^r-1$.
\end{Lem}

\begin{Proof}
If $\ord_{T-1}(P)>0$, then $P(0)=0$ and so by construction 
$$Q(T) \ =\ \frac{P(T)}{T-1}\cdot T + T\cdot(T-1)^{2^r-2}.$$
If in addition $\ord_{T-1}(P) < 2^r-1$, then $\ord_{T-1}\bigl(\frac{P(T)}{T-1}\cdot T\bigr) = \ord_{T-1}(P)-1<2^r-2$ and therefore $\ord_{T-1}(Q) = \ord_{T-1}(P)-1$. By the definition of $R(T)$ this is then also equal to $\ord_{T-1}(R)$, proving the first assertion. On the other hand, the construction of $S$ directly implies that
$$S(T) \ =\ Q(T) + T^{2^r}\cdot(Q(T) + (T-1)^{2^r-1}) 
\ =\ (T-1)^{2^r}\cdot Q(T) + T^{2^r}\cdot(T-1)^{2^r-1},$$
whence the second assertion.
\end{Proof}


\subsection{The orbit length generating function}
\label{E2}

For any integer $r\ge0$ consider the power series 
$$\Omega_r\ :=\ 
\sum_{m\ge0} \; \bigl({\textstyle\frac{X}{2}}\bigr)^m
\cdot (2Y)^{2^{m+r}-2^r} 
\ \in\ \BZ[[X,Y]].$$
For any $w\in W$ let $\Psi_w$ denote the power series from Proposition \ref{OLGFRecRels21}.

\begin{Lem}\label{bLem4}
For any polynomial $P$ in $\BF_2[T]$ of degree $< 2^r$ with $d := \ord_{T-1}(P) <2^r-1$ we have 
$$\Psi_{w_{r,P}} \ =\ 2^dY^{d+1}\Omega_r.$$
\end{Lem}

\begin{Proof}
It suffices to show the equation modulo $Y^N$ for all $N\ge0$, which we will achieve by induction on~$N$. The case $N=0$ is trivial, so assume that $N>0$ and that the equation holds universally modulo $Y^{N-1}$.

If $d>0$, then $P(1)=0$, and so $w_{r,P} = (w_{r,Q},w_{r,R})$ by Lemma \ref{bLem1}. By Proposition \ref{OLGFRecRels22} we therefore have $\Psi_{w_{r,P}} = Y\Psi_{w_{r,Q}}+Y\Psi_{w_{r,R}}$. On the other hand we have $\ord_{T-1}(Q) = \ord_{T-1}(R) = d-1$ by Lemma \ref{bLem2} and so by the induction hypothesis $\Psi_{w_{r,Q}} \equiv \Psi_{w_{r,R}} \equiv 
2^{d-1}Y^{d}\Omega_r$ modulo $Y^{N-1}$. Together this implies that $\Psi_{w_{r,P}} \equiv 2^dY^{d+1}\Omega_r$ modulo $Y^N$, as desired.

If $d=0$, then $P(1)=1$, and so $w_{r,P} = (w_{r,Q},w_{r,R})\,\sigma$ with $w_{r,Q}w_{r,R}=w_{r+1,S}$ by Lemma \ref{bLem1}. By Proposition \ref{OLGFRecRels22} we therefore have $\Psi_{w_{r,P}} = Y+XY\Psi_{w_{r+1,S}}$. Since $\ord_{T-1}(S) = 2^r-1$ by Lemma \ref{bLem2} and $2^r-1 < 2^{r+1}-1$, the induction hypothesis implies that $\Psi_{w_{r+1,S}} \equiv 2^{2^r-1}Y^{2^r}\Omega_{r+1}$ modulo $Y^{N-1}$. Together this shows that 
$$\Psi_{w_{r,P}}\ \equiv\ Y+XY2^{2^r-1}Y^{2^r}\Omega_{r+1} 
\quad\hbox{modulo\quad} Y^N.$$ 
A short calculation shows that the right hand side is equal to~$Y\Omega_r$; hence $\Psi_{w_{r,P}} \equiv Y\Omega_r$ modulo $Y^N$, as desired.
\end{Proof}

\begin{Prop}\label{bOLGF}
We have 
$$\Phi_b\ =\ 1 \ +\ \sum_{m\ge1} \; \sum_{2^{m-1}\le n<2^m} \!\!2^{n-m}X^mY^n.$$
\end{Prop}

\begin{Proof}
For $r:=1$ the polynomial $P := T$ has degree $1<2^r$ and $d := \ord_{T-1}(P)=0<2^r-1$, which satisfies the assumptions of Lemma \ref{bLem4}.
Since in this case $w_{r,P} = bb\sigma$ by (\ref{bLem0}), we find that $\Psi_{bb\sigma} = Y\Omega_1$. But by definition $b=(b,b\sigma)\sigma$, so with Proposition \ref{OLGFRecRels22} we deduce that $\Psi_b = Y+XY\Psi_{bb\sigma} = Y+XY^2\Omega_1$. A direct calculation now shows that
\UseTheoremCounterForNextEquation
\begin{equation}\label{bOLGFPsi}
\Psi_b \ =\ {\textstyle\frac{1}{2}}\cdot
\sum_{m\ge0} \; \bigl({\textstyle\frac{X}{2}}\bigr)^m
\cdot (2Y)^{2^m},
\end{equation}
and another yields the indicated formula for~$\Phi_b$.
\end{Proof}

\begin{Cor}\label{bOLGF2}
\begin{enumerate}
\item[(a)] 
The length of any orbit of $b$ on any level $n\ge0$ is the smallest power of $2$ which is greater than~$n$.
\item[(b)] For any $m\ge0$, the power $b^{2^m}$ fixes all vertices on level $2^m-1$, but none on level~$2^m$.
\end{enumerate}
\end{Cor}

\begin{Proof}
By the definition of $\Phi_b$ both assertions are equivalent to Proposition \ref{bOLGF}.
\end{Proof}


\subsection{Irrationality}
\label{E3}

\begin{Prop}\label{bNotRat}
The power series $\Phi_b$ is not a rational function of $(X,Y)$.
\end{Prop}

\begin{Proof}
By construction $\Phi_b$ is rational if and only if $\Psi_b$ is rational. If so, there exist non-zero polynomials $f,g\in\BQ[X,Y]$ with $f=g\cdot\Psi_b$. By (\ref{bOLGFPsi}) this means that
$$f(X,Y)\ =\ {\textstyle\frac{1}{2}}\cdot
\sum_{m\ge0} \; g(X,Y)\cdot\bigl({\textstyle\frac{X}{2}}\bigr)^m
\cdot (2Y)^{2^m}.$$
But for degree reasons, the summands for all $m$ with $2^m>\max\{\deg_Y(f),\deg_Y(g)\}$ cannot cancel with any other terms, yielding a contradiction. Thus $\Psi_b$ and hence $\Phi_b$ is not rational, as desired.
\end{Proof}


%
%
%
%
%
%
%
%

 \newpage


\hyphenation{polynomials}
\section{Iterated monodromy groups of quadratic polynomials: Periodic case}
\label{S=0}
\hyphenation{poly-no-mi-als}


\subsection{The iterated monodromy group}
\label{sec0IMG}

Throughout this section we fix an integer $r>0$ and a tuple $\Ux=(x_2,\ldots,x_r)$ of elements of $\{0,1\}$. Consider the elements $b_1,\ldots,b_r\in W$ defined by the recursion relations
\UseTheoremCounterForNextEquation
\begin{equation}\label{0bRecRels}
\left\{\begin{array}{ll}
b_1 = (1,b_r)\,\sigma, & \\[3pt]
b_i = (b_{i-1},1) & \hbox{for all $2\le i\le r$ with $x_i=0$,} \\[3pt]
b_i = (1,b_{i-1}) & \hbox{for all $2\le i\le r$ with $x_i=1$,} \\
\end{array}\right.
\end{equation}
and let $\GammaUx\subset W$ be the subgroup generated by them. Up to a change in notation, these are the generators and the subgroup studied by Bartholdi and Nekrashevych in \cite[\S3]{Bartholdi-Nekrashevych-2008}. Thus by 
\cite[Thm.\;5.1]{Bartholdi-Nekrashevych-2008} we have:

\begin{Thm}\label{0BarthNekrThm}
Let $f$ be any quadratic polynomial over $\BC$ and ${\eta\in\BC}$ be its unique critical point. Assume that $\eta,f(\eta),\ldots,f^{r-1}(\eta)$ are all distinct and that $f^r(\eta)=\eta$. Then the iterated monodromy group of $f$ is $W$-conjugate to~$\GammaUx$ for a certain choice of~$\Ux$.
\end{Thm}

Note that the inverses of the generators in (\ref{0bRecRels}) satisfy
$$\left\{\begin{array}{ll}
b_1^{-1} = (b_r^{-1},1)\,\sigma, & \\[3pt]
b_i^{-1} = (b_{i-1}^{-1},1) & \hbox{for all $2\le i\le r$ with $x_i=0$,} \\[3pt]
b_i^{-1} = (1,b_{i-1}^{-1}) & \hbox{for all $2\le i\le r$ with $x_i=1$.} \\
\end{array}\right.$$
Thus all the following results on $\GammaUx$ also hold if the first relation in (\ref{0bRecRels}) is replaced by the relation $b_1=(b_r,1)\,\sigma$ (see \cite[p.\;316]{Bartholdi-Nekrashevych-2008}). 
In particular, the results in the case $\Ux=(0,\ldots,0)$ apply to the subgroup generated by the elements $a_1,\ldots,a_r$ studied in \cite[\S2]{Pink2013b}, which were defined by
\UseTheoremCounterForNextEquation
\begin{equation}\label{0aRecRels}
\biggl\{\begin{array}{ll}
a_1 = (a_r,1)\,\sigma, & \\[3pt]
a_i = (a_{i-1},1) & \hbox{for all $2\le i\le r$.}\\
\end{array}
\end{equation}
Also observe:

\begin{Prop}\label{0X1-X}
The group $\Gamma_{(x_2,\ldots,x_r)}$ is conjugate to the group $\Gamma_{(1-x_2,\ldots,1-x_r)}$ under~$W$.
\end{Prop}

\begin{Proof}
Consider the element $w\in W$ that is defined by the recursion relation $w=(w,w)\,\sigma$. Then a direct calculation shows that
$$\left\{\begin{array}{ll}
wb_1^{-1}w^{-1} = (1,wb_r^{-1}w^{-1})\,\sigma, & \\[3pt]
wb_i^{-1}w^{-1} = (1,wb_{i-1}^{-1}w^{-1}) & \hbox{for all $2\le i\le r$ with $x_i=0$,} \\[3pt]
wb_i^{-1}w^{-1} = (wb_{i-1}^{-1}w^{-1},1) & \hbox{for all $2\le i\le r$ with $x_i=1$.} \\
\end{array}\right.$$
Thus the elements $w b_i^{-1}w^{-1}$ satisfy the relations (\ref{0bRecRels}) with $1-x_i$ in place of~$x_i$, and so
$w\Gamma_{(x_2,\ldots,x_r)}w^{-1} = \Gamma_{(1-x_2,\ldots,1-x_r)}$.
\end{Proof}

\bigskip
The aim of this section is to show that the orbit length generating functions of all elements of $\GammaUx$ are rational and possess an explicit common denominator. 


\subsection{Finiteness}
\label{sec0Fin}

We begin with some preparations. Let $\pi$ denote the cyclic permutation 
of the set $\{1,\ldots,r\}$ defined by
\UseTheoremCounterForNextEquation
\begin{equation}\label{0PiDef}
\pi(i)\ :=\ 
\biggl\{\begin{array}{ll}
r & \hbox{if $i=1$,}\\[3pt]
i-1 & \hbox{if $i>1$.}
\end{array}
\end{equation}
Then the recursion relations (\ref{0bRecRels}) express each $b_i$ in terms of $b_{\pi(i)}$.

\begin{Def}\label{0LengthDef}
The \emph{length} $|w|$ of an element $w\in\GammaUx$ is the minimal length of a word over the alphabet $\{b_1^{\pm1},\ldots,b_r^{\pm1}\}$ that represents~$w$. Any word of minimal length representing $w$ is called a \emph{minimal word for~$w$}.
\end{Def}

\begin{Lem}\label{0Lem1}
For any element $w=(u,v)\,\sigma^\mu \in\GammaUx$ we have $u,v\in\GammaUx$ and 
$${|uv| \le |u|+|v| \le |w|}.$$
\end{Lem}

\begin{Proof}
By the recursion relations (\ref{0bRecRels}), any letter $b_i^{\pm1}$ in a minimal word for $w$ contributes precisely one letter $b_{\pi(i)}^{\pm1}$ to a word representing $u$ or~$v$.
This implies the second inequality, and the first one follows directly from the definition of length.
\end{Proof}

\begin{Lem}\label{0Lem2}
For all $w\in\GammaUx$ and all $w'\in\Desc(w)$ we have $w'\in\GammaUx$ with $|w'|\le|w|$.
\end{Lem}

\begin{Proof}
By Definition \ref{DescDef} and iteration this follows from Lemma \ref{0LengthDef}.
\end{Proof}

\begin{Prop}\label{0OLGFFin}
Every element of $\GammaUx$ is \OLGF-finite.
\end{Prop}

\begin{Proof}
Since $\GammaUx$ contains only finitely many elements of any given length, Lemma \ref{0Lem2} implies that $\Desc(w)$ is finite for any $w\in\GammaUx$, as desired.
\end{Proof}

\medskip
Combining Proposition \ref{0OLGFFin} with Theorem \ref{OLGFRat} we find that the orbit length generating functions of all elements of $\GammaUx$ are rational. By Theorem \ref{OLGFLinComb} the study of their denominators reduces to the case of \OLGF-irreducible elements. 


\subsection{Properties of \OLGF-irreducible elements}
\label{sec0Prop}

\begin{Lem}\label{0Lem3}
Any \OLGF-irreducible element $w\in\GammaUx$ has a unique first descendant $w'$ which is \OLGF-irreducible with $|w'|=|w|$. Moreover $w$ is either
$W$-conjugate to $(w',1)\,\sigma$, or equal to $(w',1)$ or $(1,w')$.
\end{Lem}

\begin{Proof}
Suppose first that $w=(u,v)\,\sigma$. Then $w$ is $W$-conjugate to $(uv,1)\,\sigma$, and $uv$ is the unique first descendant of~$w$. Thus the assumption $w\in\Desc(w)$ means that $w$ is equal to or a descendant of~$uv$. On the one hand this implies that $uv$ is a descendant of itself; hence $uv$ is \OLGF-irreducible. On the other hand it implies by Lemma \ref{0Lem2} that $|w|\le|uv|\le|w|$ and hence $|uv|=|w|$, and we are done with $w':=uv$.

Suppose now that $w=(u,v)$, so that $u$ and $v$ are the first descendants of~$w$. Then the assumption $w\in\Desc(w)$ means that $w$ is equal to, or a descendant of, one of $u$, $v$; let us call it~$w'$. On the one hand this implies that $w'$ is a descendant of itself; hence $w'$ is \OLGF-irreducible. On the other hand it implies by Lemma \ref{0Lem2} that $|w|\le|w'|\le|w|$ and hence $|w'|=|w|$. Plugging this into the inequality $|u|+|v|\le|w|$ from Lemma \ref{0Lem1}, we now deduce that the other entry of $(u,v)$ has length $0$ and is therefore the identity element. Thus $w=(w',1)$ or $w=(1,w')$. This makes $w'$ unique (though for $w'=1$ we can write $w$ in both ways). Since $(1,w')$ is $W$-conjugate to $(w',1)$, in either case we are done.
\end{Proof}

\medskip
Next we look at signs. The same proof as that of \cite[Prop.\;2.1.1]{Pink2013b} shows:

\begin{Lem}\label{0Lem4}
For all $n\ge1$ and all $1\le i\le r$ we have
$$\sgn_n(b_i) \ =\ 
\biggl\{\!\begin{array}{rl}
-1 & \hbox{if $n    \equiv i\bmod r$,}\\[3pt]
 1 & \hbox{if $n\not\equiv i\bmod r$.}
\end{array}$$
Thus for any fixed $w\in\GammaUx$, the value $\sgn_n(w)$ for $n\ge1$ depends only on $n\bmod r$.
\end{Lem}

To any element $w\in\GammaUx$ we associate the subset
\UseTheoremCounterForNextEquation
\begin{equation}\label{0JwDef}
J_w \ :=\ \{1\le i\le r\mid\sgn_i(w)=-1\}.
\end{equation}

\begin{Lem}\label{0JwInd}
For any $w$ and $w'$ as in Lemma \ref{0Lem3} we have $J_{w'}=\pi(J_w)$.
\end{Lem}

\begin{Proof}
The recursion relations for signs and their invariance under conjugation implies that $\sgn_i(w) = \sgn_{i-1}(w')$ for all $i\ge2$. Using the periodicity from Lemma \ref{0Lem4} we also find that $\sgn_1(w) = \sgn_{r+1}(w) = \sgn_r(w')$. By (\ref{0PiDef}) we therefore have $\sgn_i(w) = \sgn_{\pi(i)}(w')$ for all $1\le i\le r$, or equivalently $J_{w'}=\pi(J_w)$.
\end{Proof}


\subsection{Conjugacy classes of \OLGF-irreducible elements}
\label{sec0Conj}

\begin{Lem}\label{0Lem19}
Consider any distinct indices $i_1,\ldots,i_k \in \{1,\ldots,r\}$, in any order. Set $\mu:=1$ if $1$ appears among them, and $\mu:=0$ otherwise. Then $a_{i_1}\cdots a_{i_k}$ is conjugate to $(a_{\pi(i_1)}\cdots a_{\pi(i_k)},1)\,\sigma^\mu$ under~$W$.
\end{Lem}

\begin{Proof}
If $1$ does not appear among $i_1,\ldots,i_k$, the recursion relations (\ref{0aRecRels}) imply that 
$$a_{i_1}\cdots a_{i_k} \ =\  
(a_{i_1-1},1)\cdots(a_{i_k-1},1) \ =\ 
(a_{\pi(i_1)}\cdots a_{\pi(i_k)},1)\,\sigma^\mu,$$
and the assertion follows. Otherwise let $j$ be the unique index with $i_j=1$. Then the recursion relations (\ref{0aRecRels}) imply that 
$$\begin{array}{rl}
a_{i_1}\cdots a_{i_k} \ =&
(a_{i_1-1},1)\cdots(a_{i_{j-1}-1},1)
\cdot(a_r,1)\,\sigma \cdot
(a_{i_{j+1}-1},1)\cdots(a_{i_k-1},1) \\[3pt]
\ =& (a_{\pi(i_1)}\cdots a_{\pi(i_j)},
a_{\pi(i_{j+1})}\cdots a_{\pi(i_k)})\,\sigma^\mu.
\end{array}$$
This is $W$-conjugate to $(a_{\pi(i_1)}\cdots\cdots a_{\pi(i_k)},1)\,\sigma^\mu$, as desired.
\end{Proof}


\begin{Prop}\label{0IrredConj}
Consider any \OLGF-irreducible element $w\in\GammaUx$. Let $i_1,\ldots,i_k$ be the distinct elements of~$J_w$, in any order. Then $w$ is conjugate to $a_{i_1}\cdots a_{i_k}$ under~$W$.
\end{Prop}

\begin{Proof}
By \cite[Lemma\;1.3.3]{Pink2013b} it suffices to prove that the restrictions $w|_{T_n}$ and $a_{i_1}\cdots a_{i_k}|_{T_n}$ are conjugate in the automorphism group of $T_n$ for every $n\ge0$. We will achieve this by induction on~$n$. For $n=0$ the assertion is trivially true, so assume that $n>0$ and that the assertion is universally true for the restrictions to $T_{n-1}$.

Let $w'\in\GammaUx$ be the unique \OLGF-irreducible descendant of~$w$ from Lemma \ref{0Lem3}. Then $w$ is conjugate to $(w',1)\,\sigma^\mu$ for some $\mu\in\{0,1\}$. Thus $\sgn_1(w)=(-1)^\mu$, and hence $\mu=1$ if and only if $1\in J_w$. Also, Lemma \ref{0JwInd} shows that $\pi(i_1),\ldots,\pi(i_k)$ are the distinct elements of~$J_{w'}$.
By the induction hypothesis $w'|_{T_{n-1}}$ is therefore conjugate to $a_{\pi(i_1)}\cdots a_{\pi(i_k)}|_{T_{n-1}}$ under the automorphism group of~$T_{n-1}$. Thus $w|T_n$ is conjugate to $(a_{\pi(i_1)}\cdots a_{\pi(i_k)},1)\,\sigma^\mu|_{T_n}$ under the automorphism group of~$T_n$. From Lemma \ref{0Lem19} it now follows that $w|T_n$ is conjugate to $a_{i_1}\cdots a_{i_k}|_{T_n}$ under the automorphism group of~$T_n$, as desired.
\end{Proof}


\medskip
The next result concerns the following condition on a subset $J\subset\{1,\ldots,r\}$.

\begin{Cond}\label{0JCond}
For any $n\ge0$ with $1\not\in\pi^n(J)$, the values $x_i$ for all $i\in \pi^n(J)$ are equal.
\end{Cond}

\begin{Prop}\label{0IrredExists}
For any subset $J\subset\{1,\ldots,r\}$ satisfying Condition \ref{0JCond} there exists a \OLGF-irreducible element $w\in\GammaUx$ with $J_w=J$.
\end{Prop}

\begin{Proof}
Consider any integer $n\ge0$. For the purpose of this proof we call any element of $\GammaUx$ of the form $b_{i_1}\cdots b_{i_k}$, where $i_1,\ldots,i_k$ are the distinct elements of $\pi^n(J)$ in any order, \emph{strongly of type $\pi^n(J)$}. We claim that any element that is strongly of type $\pi^n(J)$ possesses a first descendant which is strongly of type $\pi^{n+1}(J)$.

Granting this, by induction on $n$ it follows that for any $n\ge1$, any element that is strongly of type $J$ possesses a descendant which is strongly of type $\pi^n(J)$. Since $\pi$ is a permutation of finite order, we deduce that any element that is strongly of type $J$ possesses a descendant which is again strongly of type~$J$. As there are only finitely many elements that are strongly of type $J$, and being a descendant is a transitive relation, it follows that some element $w_0$ that is strongly of type $J$ must be its own descendant. This element is therefore \OLGF-irreducible. Finally, writing $w_0=b_{i_1}\cdots b_{i_k}$ where $i_1,\ldots,i_k$ are the distinct elements of~$J$, Lemma \ref{0Lem4} implies that $J_w=J$, as desired.

To prove the claim consider $w := b_{i_1}\cdots b_{i_k}$ where $i_1,\ldots,i_k$ are the distinct elements of $\pi^n(J)$. 
Suppose first that $1\not\in\pi^n(J)$. Then by Condition \ref{0JCond} the values $x_i$ are equal for all $i\in\pi^n(J)$. Thus the recursion relations (\ref{0bRecRels}) imply that $b_{i_1}\cdots b_{i_k} = (b_{\pi(i_1)}\cdots b_{\pi(i_k)},1)$ or $(1,b_{\pi(i_1)}\cdots b_{\pi(i_k)})$. In both cases $w$ has the first descendant $b_{\pi(i_1)}\cdots b_{\pi(i_k)}$, which is strongly of type $\pi^{n+1}(J)$.

Suppose now that $1\in\pi^n(J)$. Then $\sgn_1(w)=-1$ by Lemma \ref{1LemSigns} (a); hence $w$ has the form $w=(u,v)\,\sigma$. By the recursion relations (\ref{0bRecRels}), any factor $b_{i_j}$ of $w=b_{i_1}\cdots b_{i_k}$ contributes precisely one factor $b_{\pi(i_j)}$ to the product~$uv$. Thus $uv$ is a product of the elements $b_{\pi(i_1)},\ldots,b_{\pi(i_k)}$ in some order. It is therefore strongly of type $\pi^{n+1}(J)$, as desired.
\end{Proof}


\begin{Prop}\label{0IrredConjToSplit}
Any \OLGF-irreducible element $w$ of $\GammaUx$ is $W$-conjugate to a \OLGF-irreducible element of $\Gamma_{(0,\ldots,0)}$.
\end{Prop}

\begin{Proof}
By Proposition \ref{0IrredConj} it is conjugate to $a_{i_1}\cdots a_{i_k} \in \Gamma_{(0,\ldots,0)}$, where $i_1,\ldots,i_k$ are the distinct elements of~$J_w$ in any order. But the same argument as in the proof of Proposition \ref{0IrredExists} shows that for some order, the element $a_{i_1}\cdots a_{i_k}$ is \OLGF-irreducible.
\end{Proof}


\subsection{Some combinatorics}
\label{sec0Comb}

The content of this subsection and the next is needed only to determine the precise lowest common denominator in Theorem \ref{0MainThm} below, and can be skipped if one is happy with some common denominator.

\medskip
For all $x\in\{0,1\}$ we set 
\UseTheoremCounterForNextEquation
\begin{eqnarray}
\label{1CombSxDef}
S^x &\!\!:=\!\!& \{2\le i\le r\mid x_i=x\} \quad\hbox{and} \\[3pt]
\UseTheoremCounterForNextEquation
\label{1CombIrxDef}
I_r^x &\!\!:=\!\!& \pi(S^x).
\end{eqnarray}
For all  $x\in\{0,1\}$ and $1<i\le r$ we define by descending induction
\UseTheoremCounterForNextEquation
\begin{equation}\label{1CombIixDef}
I_{i-1}^x\ :=\ 
\biggl\{\begin{array}{ll}
\pi(I_i^x \cap S^{x_i}) & \hbox{if $1\not\in I_i^x$,}\\[5pt]
\pi(I_i^x \cup S^{1-x_i}) & \hbox{if $1\in I_i^x$.}
\end{array}
\end{equation}

\begin{Lem}\label{1CombLem1}
For all $1\le i\le r$ we have a decomposition into disjoint subsets
$$\{1,\ldots,r\}\ =\ \{i\} \sqcup I_i^0 \sqcup I_i^1.$$
\end{Lem}

\begin{Proof}
{}From (\ref{1CombSxDef}) we deduce that
$\{1,\ldots,r\} = \{1\} \sqcup S^0 \sqcup S^1$.
By (\ref{1CombIrxDef}) this implies the desired assertion for $i=r$. Suppose that the assertion holds for $1<i\le r$.
Then there is a unique index $x\in\{0,1\}$ with $1\in I_i^x$ and $1\not\in I_i^{1-x}$. By (\ref{1CombIixDef}) we thus have
$$\begin{array}{ll}
I_{i-1}^{1-x} \!\!\! &=\; \pi(I_i^{1-x} \cap S^{x_i}) 
\quad\hbox{and}\\[5pt]
I_{i-1}^x     \!\!\! &=\; \pi(I_i^x \cup S^{1-x_i}).
\end{array}$$
The fact that $1\not\in I_i^{1-x}$ also implies that
$$I_i^{1-x} \ =\ (I_i^{1-x} \cap S^{1-x_i})
\sqcup (I_i^{1-x} \cap S^{x_i}).$$
The induction hypothesis and the fact that $i\not\in S^{1-x_i}$ imply that
$$I_i^x \cup (I_i^{1-x} \cap S^{1-x_i})
\ =\ I_i^x \cup S^{1-x_i}.$$
Together it follows that
\begin{eqnarray*}
\{1,\ldots,r\} 
&\!\!=\!\!& \pi\bigl(\{i\} \sqcup I_i^x \sqcup I_i^{1-x} \bigr) \\[3pt]
&\!\!=\!\!& \pi\bigl(\{i\} \sqcup I_i^x \sqcup (I_i^{1-x} \cap S^{1-x_{i}}) \sqcup (I_i^{1-x} \cap S^{x_i}) \bigr) \\[3pt]
&\!\!=\!\!& \pi\bigl(\{i\} \sqcup (I_i^x \cup S^{1-x_i}) \sqcup (I_i^{1-x} \cap S^{x_i}) \bigr) \\[3pt]
&\!\!=\!\!& \{i-1\} \sqcup I_{i-1}^x \sqcup I_{i-1}^{1-x},
\end{eqnarray*}
and the desired assertion holds for~$i-1$. By downward induction it follows for all~$i$.
\end{Proof}

\begin{Lem}\label{1CombLem2}
For any distinct $1\le i,j\le r$ there exist $x,y\in\{0,1\}$ such that 
$$I_i^x \cup I_j^y \ =\ \{1,\ldots,r\}.$$
\end{Lem}

\begin{Proof}
Suppose first that one of $i$, $j$ is equal to~$r$. By symmetry we may assume that $i<j=r$. By Lemma \ref{1CombLem1} there is a unique $x\in\{0,1\}$ such that $1\in I_{i+1}^x$. With $y := x_{i+1}$ the constructions (\ref{1CombIrxDef}) and (\ref{1CombIixDef}) then imply that
$$I_i^x \cup I_r^y \ =\ \pi\bigl(I_{i+1}^x \cup S^{1-y} \cup S^y\bigr).$$
Since $1\in I_{i+1}^x$ and $\{1\} \cup S^{1-y} \cup S^y = \{1,\ldots,r\}$, the right hand side is equal to $\{1,\ldots,r\}$, as desired.

Suppose now that the assertion holds for given $i,j>1$. We then prove it for $i-1$ and $j-1$.  By Lemma \ref{1CombLem1} there are unique $x,y\in\{0,1\}$ such that $1\in I_i^x\cap I_j^y$. The construction (\ref{1CombIixDef}) then implies that
\UseTheoremCounterForNextEquation
\begin{equation}\label{1CombLem2Disp}
I_{i-1}^x \cup I_{j-1}^y \ =\ 
\pi\bigl(I_i^x \cup S^{1-x_i} \cup I_j^y \cup S^{1-x_j} \bigr).
\end{equation}
If $x_i\not=x_j$, the right hand side of (\ref{1CombLem2Disp}) contains $\pi\bigl(\{1\}\cup S^0 \cup S^1 \bigr) = \{1,\ldots,r\}$, and we are done. Otherwise abbreviate $z := x_i=x_j$. Using the induction hypothesis choose $x',y'\in\{0,1\}$ such that $I_i^{x'} \cup I_j^{y'} = \{1,\ldots,r\}$. Then in particular $1\in I_i^{x'} \cup I_j^{y'}$, and so either $x'=x$ or $y'=y$ or both. If $(x',y')=(x,y)$, the right hand side of (\ref{1CombLem2Disp}) contains $\pi\bigl(I_i^x \cup I_j^y\bigr) = \{1,\ldots,r\}$, and we are done. Otherwise by symmetry we may without loss of generality assume that $(x',y')=(x,1-y)$. Instead of $I_{i-1}^x \cup I_{j-1}^y$ we then look at $I_{i-1}^x \cup I_{j-1}^{1-y}$. Since $1\in I_i^x\setminus I_j^{1-y}$, the construction (\ref{1CombIixDef}) implies that
$$I_{i-1}^x \cup I_{j-1}^{1-y} \ =\ 
\pi\bigl(I_i^x \cup S^{1-z} \cup (I_j^{1-y} \cap S^z) \bigr).$$
Since $1\in I_i^x$ and $\{1\} \cup S^{1-z} \cup S^z = \{1,\ldots,r\}$, we deduce that
$$I_{i-1}^x \cup I_{j-1}^{1-y} \ \supset\ \pi\bigl(I_i^x \cup I_j^{1-y}\bigr) \ =\ \{1,\ldots,r\},$$
and again we are done. The lemma thus follows by descending induction.
\end{Proof}

\begin{Lem}\label{1CombLem3}
There exist $1\le i\le r$ and $x\in\{0,1\}$ such that $I_i^x=\emptyset$.
\end{Lem}

\begin{Proof}
Choose $i$ and $x$ such that $|I_i^x|$ is minimal. If $|I_i^x|>0$, pick any $j\in I_i^x$. Then Lemma \ref{1CombLem1} implies that $j\not=i$. 
Using Lemma \ref{1CombLem2} choose $x',y\in\{0,1\}$ such that 
$I_i^{x'} \cup I_j^y = \{1,\ldots,r\}$. 
Then by Lemma \ref{1CombLem1} for $j$ in place of $i$ we have $j\not\in I_j^y$, and therefore $j\in I_i^{x'}$. Thus $j\in I_i^{x'} \cap I_i^x$, which by Lemma \ref{1CombLem1} implies that $x'=x$.
Therefore $I_i^x \cup I_j^y = \{1,\ldots,r\}$. 
Counting elements, and using Lemma \ref{1CombLem1} for $j$ in place of $i$ again, we deduce that
$$|I_i^x| + |I_j^y| \ \ge\ r \ =\ 1 + |I_j^{1-y}| + |I_j^y|.$$
Therefore $|I_i^x| \ge 1 + |I_j^{1-y}| > |I_j^{1-y}|$, contradicting the minimality of~$|I_i^x|$. Thus after all we have $|I_i^x|=0$, and hence $I_i^x=\emptyset$, as desired.
\end{Proof}

\begin{Lem}\label{1CombLem4}
There exists $x\in\{0,1\}$ such that $I_1^x=\emptyset$.
\end{Lem}

\begin{Proof}
By Lemma \ref{1CombLem3} there exists a smallest index $1\le i\le r$ such that $I_i^x=\emptyset$ for some $x\in\{0,1\}$. If that index is $>1$, we in particular have $1\not\in I_i^x$; hence the construction (\ref{1CombIixDef}) implies that $I_{i-1}^x = \pi(I_i^x \cap S^{x_i}) = \emptyset$, contradicting the minimality of~$i$. 
\end{Proof}


\subsection{Minimal words for \OLGF-irreducible elements}
\label{sec0More}

In this subsection we study the minimal words for \OLGF-irreducible elements in more detail. 

\medskip
Here and only here we use the following abbreviations: For any subset $I\subset\{1,\ldots,r\}$ we let $\langle I\rangle$ denote any (possibly empty) word over the alphabet $\{b_i^{\pm1}\mid i\in I\}$. A concatenation of expressions $\langle I\rangle$ for subsets $I$ and/or of individual letters $b_i^{\pm1}$ represents the concatenation of any words or letters of the indicated form. An overline $\overline{\phantom{iiiiiiiii}}$ over such a 
confounded expression means that the template is repeated an arbitrary non-negative number of times.
One should keep in mind that this notation refers to words and not to the group elements represented by them.

\begin{Lem}\label{0LemIJ}
If $w\in\GammaUx$ is represented by a word of the form $\langle I\rangle$, then $J_w\subset I$.
\end{Lem}

\begin{Proof}
Lemma \ref{0Lem4} implies that $\sgn_i(b_i)=-1$ and ${\sgn_i(b_j)=1}$ whenever $i\not=j$.
\end{Proof}

\medskip
Let $S^x$ and $I_i^x$ be as in the preceding subsection. 

\begin{Lem}\label{0MoreLem1}
Consider any \OLGF-irreducible element $w\in\GammaUx$ with $1\not\in J_w$. Then any minimal word for $w$ has one of the forms
$$\begin{array}{c}
\langle S^0\rangle\ \overline{b_1^{-1} \langle S^1\rangle\,b_1 \langle S^0\rangle}, \\[5pt]
\langle S^1\rangle\ 
\overline{b_1 \langle S^0\rangle\,b_1^{-1} \langle S^1\rangle}.
\end{array}$$
\end{Lem}

\begin{Proof}
The assumption $1\not\in J_w$ means that $\sgn_1(w)=1$. Thus $w=(u,v)$ for certain $u,v\in\GammaUx$. By the recursion relations (\ref{0bRecRels}), any letter $b_k^{\pm1}$ in the minimal word for $w$ contributes precisely one letter $b_{\pi(k)}^{\pm1}$ to a word representing $u$ or~$v$. Since one of $u$, $v$ has the same length as $w$ by Lemma \ref{0Lem3}, this letter must always land in the same one of $u$,~$v$. 
Now suppose that the minimal word in question contains a subword of one of the forms
\smallskip
$$\begin{array}{|c|c|c|c|}
\hline
{\Large\strut} 
  b_1b_1
& b_1^{-1}b_i^{\pm1}
& b_1b_j^{\pm1} 
& b_i^{\pm1}b_j^{\pm1} \\[3pt]
\hline
{\Large\strut} 
  b_1^{-1}b_1^{-1}
& b_i^{\pm1}b_1
& b_j^{\pm1}b_1^{-1}
& b_j^{\pm1}b_i^{\pm1} \\[3pt]
\hline
\end{array}$$

\medskip\noindent
for some $i\in S^0$ and $j\in S^1$ and independent exponents $\pm1$. By (\ref{0bRecRels}) the recursive expansion of this subword is, respectively:
\smallskip
$$\begin{array}{|c|c|c|c|}
\hline
{\Large\strut} 
  (b_r,b_r)
& (b_r^{-1},b_{i-1}^{\pm1})\,\sigma
& (b_{j-1}^{\pm1},b_r)\,\sigma
& (b_{i-1}^{\pm1},b_{j-1}^{\pm1}) \\[3pt]
\hline
{\Large\strut} 
  (b_r^{-1},b_r^{-1})
& (b_{i-1}^{\pm1},b_r)\,\sigma
& (b_r^{-1},b_{j-1}^{\pm1})\,\sigma
& (b_{i-1}^{\pm1},b_{j-1}^{\pm1}) \\[3pt]
\hline
\end{array}$$

\medskip\noindent
Thus the two letters of this subword bequeath one letter to each of $u$ and~$v$, yielding a contradiction. Therefore the minimal word does not contain a subword of the above form.
This means that the minimal word is a subword of a word of the form
$$\ldots b_1 \langle S^0\rangle\, 
    b_1^{-1} \langle S^1\rangle\, 
         b_1 \langle S^0\rangle\, 
    b_1^{-1} \langle S^1\rangle\, b_1 \ldots.$$
Finally, since $\sgn_1(w)=1$, and $\sgn_1(b_i)=-1$ only for $i=1$, the total number of letters $b_1^{\pm1}$ is even. Depending on the first letter the minimal word therefore has the indicated form.
\end{Proof}

\medskip
For the following argument we fix a \OLGF-irreducible element $w_0\in\GammaUx$ with $1\not\in J_w$. We construct a sequence of \OLGF-irreducible elements $w_n\in\GammaUx$ by defining each $w_{n+1}$ as the first descendant of $w_n$ furnished by Lemma \ref{0Lem3}. 
We also fix any minimal word $\tildew_0$ for~$w_0$. By repeated recursive expansion using the relations (\ref{0bRecRels}) this yields a minimal word $\tildew_n$ for $w_n$ for every $n\ge0$.

\begin{Lem}\label{0JIter}
For every $n\ge0$ we have $J_{w_n} = \pi^n(J_{w_0})$.
\end{Lem}

\begin{Proof}
This follows by induction from Lemma \ref{0JwInd}.
\end{Proof}

\begin{Lem}\label{0MoreLem2}
For any $1\le i\le r$ the word $\tildew_{r+1-i}$ has one of the forms
$$\begin{array}{c}
\langle I_i^0\rangle\ \overline{b_i^{-1} \langle I_i^1\rangle\,b_i \langle I_i^0\rangle}, \\[5pt]
\langle I_i^1\rangle\ 
\overline{b_i \langle I_i^0\rangle\,b_i^{-1} \langle I_i^1\rangle}.
\end{array}$$
\end{Lem}

\begin{Proof}
Recall from (\ref{1CombIrxDef}) that $\pi(S^x)=I_r^x$ for each $x=0,1$. Thus the recursion relations (\ref{0bRecRels}) show that any word of the form $\langle S^0\rangle\ \overline{b_1^{-1} \langle S^1\rangle\,b_1 \langle S^0\rangle}$ expands to one of the form $\bigl(\langle I_r^0\rangle\ \overline{b_r^{-1} \langle I_r^1\rangle\,b_r \langle I_r^0\rangle},\ 1\bigr)$, and any word of the form $\langle S^1\rangle\ \overline{b_1 \langle S^0\rangle\,b_1^{-1} \langle S^1\rangle}$ expands to one of the form $\bigl(1,\ \langle I_r^1\rangle\ \overline{b_r \langle I_r^0\rangle\,b_r^{-1} \langle I_r^1\rangle}\,\bigr)$. Lemma \ref{0MoreLem1} therefore implies the desired assertion in the case $i=r$ for the word $\tildew_{r+1-i}=\tildew_1$.

Suppose now that the assertion holds for some $1<i\le r$. We then prove it for $i-1$. We first look at the individual pieces of $\tildew_{r+1-i}$. 
Using Lemma \ref{1CombLem1} let $x\in\{0,1\}$ be the unique index with $1\in I_i^{1-x}$ and $1\not\in I_i^x$.

\begin{Sublem}\label{0MoreLem2a}
The recursive expansion of the letter $b_i$ is
$$\begin{array}{ll}
(b_{i-1},1) & \hbox{if $x_i=0$,} \\[3pt]
(1,b_{i-1}) & \hbox{if $x_i=1$.} \\
\end{array}$$
The recursive expansion of any word of the form $\langle I_i^x\rangle$ has the form
$$\begin{array}{ll}
\bigl(\langle I_{i-1}^x\rangle, \langle I_{i-1}^{1-x}\rangle \bigr)
& \hbox{if $x_i=0$,} \\[3pt]
\bigl(\langle I_{i-1}^{1-x}\rangle, \langle I_{i-1}^x\rangle \bigr)
& \hbox{if $x_i=1$.}
\end{array}$$
The recursive expansion of any word of the form $\langle I_i^{1-x}\rangle$ has one of the forms
$$\begin{array}{l}
\bigl(\langle I_{i-1}^{1-x}\rangle, \langle I_{i-1}^{1-x}\rangle \bigr), \\[5pt]
\bigl(\langle I_{i-1}^{1-x}\rangle, \langle I_{i-1}^{1-x}\rangle \bigr)
\,\sigma.
\end{array}$$
\end{Sublem}

\begin{Proof}
The first statement is a special case of the recursion relations (\ref{0bRecRels}). 
Since $1\not\in I_i^x$, the relations also imply that the recursive expansion of any word of the form $\langle I_i^x\rangle$ has the form $\bigl(\langle\pi(I_i^x\cap S^0)\rangle, \langle\pi(I_i^x\cap S^1)\rangle \bigr)$. But by (\ref{1CombIixDef}) we have $\pi(I_i^x \cap S^{x_i})= I_{i-1}^x$ and $\pi(I_i^x \cap S^{1-x_i}) \subset I_{i-1}^{1-x}$, so the second statement follows. Likewise the recursive expansion of any word of the form $\langle I_i^{1-x}\rangle$ involves only letters $b_{\pi(j)}^{\pm1}$ for $j\in I_i^{1-x}$ and (possibly) some factors~$\sigma$. It is therefore of the form $\bigl(\langle\pi(I_i^{1-x})\rangle, \langle\pi(I_i^{1-x})\rangle \bigr)$ or $\bigl(\langle\pi(I_i^{1-x})\rangle, \langle\pi(I_i^{1-x})\rangle \bigr)\,\sigma$. Since $1\in I_i^{1-x}$, by (\ref{1CombIixDef}) we have $\pi(I_i^{1-x}) \subset I_{i-1}^{1-x}$, and the third statement follows.
\end{Proof}

\medskip
Returning to the proof of Lemma \ref{0MoreLem2}, we now set $\mu := 1-2x$. Then by the induction hypothesis $\tildew_{r+1-i}$ has one of the forms
\UseTheoremCounterForNextEquation
\begin{eqnarray}
\label{0MoreLem2Form1}
& \langle I_i^{1-x}\rangle\ \overline{b_i^\mu \langle I_i^x\rangle\,b_i^{-\mu} \cdot \langle I_i^{1-x}\rangle}, & \\
\UseTheoremCounterForNextEquation
\label{0MoreLem2Form2}
& \langle I_i^x\rangle\ \overline{b_i^{-\mu} \langle I_i^{1-x}\rangle\,b_i^\mu \cdot \langle I_i^x\rangle}, &
\end{eqnarray}
and we must prove the same for $\tildew_{r+2-i}$ with $i-1$ in place of~$i$.

In the case (\ref{0MoreLem2Form1}) Sublemma \ref{0MoreLem2a} implies that the recursive expansion of $\tildew_{r+1-i}$ is a product of terms of the form $\bigl(b_{i-1}^\mu\langle I_{i-1}^x\rangle\,b_{i-1}^{-\mu}, \langle I_{i-1}^{1-x}\rangle \bigr)$ 
or $\bigl(\langle I_{i-1}^{1-x}\rangle, b_{i-1}^\mu\langle I_{i-1}^x\rangle\,b_{i-1}^{-\mu} \bigr)$ 
and/or $\bigl(\langle I_{i-1}^{1-x}\rangle, \langle I_{i-1}^{1-x}\rangle \bigr)$ 
and/or $\sigma$. It is thus equal to $(\tilde u,\tilde v)$ or $(\tilde u,\tilde v)\,\sigma$, where both $\tilde u$ and $\tilde v$ are products of terms of the form $b_{i-1}^\mu\langle I_{i-1}^x\rangle\,b_{i-1}^{-\mu}$
and/or $\langle I_{i-1}^{1-x}\rangle$. The next descendant $\tildew_{r+2-i}$ is equal to $\tilde u$ or $\tilde v$ or~$\tilde u\tilde v$ and therefore also such a product. Thus the lemma follows for $i-1$. 

In the case (\ref{0MoreLem2Form2}) suppose first that $\tildew_{r+1-i}$ does not contain the letter $b_i^{\pm1}$. Then it has the form 
$\langle I_i^x\rangle$. By Sublemma \ref{0MoreLem2a} it thus has the recursive expansion $(\tilde u,\tilde v)$ with one entry of the form $\langle I_{i-1}^x\rangle$ and the other of the form $\langle I_{i-1}^{1-x}\rangle$. Since in this case the next descendant $\tildew_{r+2-i}$ is equal to one of $\tilde u$, $\tilde v$, the lemma again follows for $i-1$. 

Now suppose that $\tildew_{r+1-i}$ has the form (\ref{0MoreLem2Form2}) and contains the letter $b_i^{\pm1}$. We then regroup its factors in the form 
\UseTheoremCounterForNextEquation
\begin{equation}\label{0MoreLem2Form3}
\underbrace{{\large\strut}\langle I_i^x\rangle\,b_i^{-\mu}} \cdot 
\underbrace{{\large\strut}\langle I_i^{1-x}\rangle\ \overline{b_i^\mu \langle I_i^x\rangle\,b_i^{-\mu} \cdot \langle I_i^{1-x}\rangle}\,}
\cdot\, \underbrace{{\large\strut}b_i^\mu \langle I_i^x\rangle}.
\end{equation}
As in the case (\ref{0MoreLem2Form1}) the whole shebang in the middle expands to $(\tilde u,\tilde v)$ or $(\tilde u,\tilde v)\,\sigma$, where $\tilde u$ and $\tilde v$ are products of terms of the form $b_{i-1}^\mu\langle I_{i-1}^x\rangle\,b_{i-1}^{-\mu}$ and/or $\langle I_{i-1}^{1-x}\rangle$. 
Assume first that it expands to $(\tilde u,\tilde v)$. Sublemma \ref{0MoreLem2a} then implies that $\tildew_{r+1-i}$ expands to 
$$\begin{array}{ll}
\bigl( \langle I_{i-1}^x\rangle\,b_{i-1}^{-\mu}
\cdot\tilde u\cdot 
b_{i-1}^\mu\langle I_{i-1}^x\rangle,\ 
\langle I_{i-1}^{1-x}\rangle 
\cdot\tilde v\cdot
\langle I_{i-1}^{1-x}\rangle \bigr)
& \hbox{if $x_i=0$,} \\[5pt]
\bigl( \langle I_{i-1}^{1-x}\rangle 
\cdot\tilde u\cdot
\langle I_{i-1}^{1-x}\rangle,\ 
\langle I_{i-1}^x\rangle\,b_{i-1}^{-\mu}
\cdot\tilde v\cdot 
b_{i-1}^\mu\langle I_{i-1}^x\rangle \bigr)
& \hbox{if $x_i=1$.}
\end{array}$$
By construction with Lemma \ref{0Lem3} the next descendant $\tildew_{r+2-i}$ is the unique non-empty entry of this pair; hence it is the one containing $b_{i-1}^{\pm1}$. Since $\tilde u$ and $\tilde v$ are products of terms of the form $b_{i-1}^\mu\langle I_{i-1}^x\rangle\,b_{i-1}^{-\mu}$ and/or $\langle I_{i-1}^{1-x}\rangle$, both
$\langle I_{i-1}^x\rangle\,b_{i-1}^{-\mu} \cdot\tilde u\cdot b_{i-1}^\mu\langle I_{i-1}^x\rangle$
and $\langle I_{i-1}^x\rangle\,b_{i-1}^{-\mu} \cdot\tilde v\cdot b_{i-1}^\mu\langle I_{i-1}^x\rangle$
are products of terms of the form $\langle I_{i-1}^x\rangle$ and/or $b_{i-1}^{-\mu}\langle I_{i-1}^{1-x}\rangle\,b_{i-1}^\mu$. Thus $\tildew_{r+2-i}$ is such a product, and so the lemma holds for $i-1$. 

Assume now that the middle of (\ref{0MoreLem2Form3}) expands to $(\tilde u,\tilde v)\,\sigma$. Sublemma \ref{0MoreLem2a} then implies that $\tildew_{r+1-i}$ expands to 
$$\begin{array}{ll}
\bigl( \langle I_{i-1}^x\rangle\,b_{i-1}^{-\mu}
\cdot\tilde u\cdot 
\langle I_{i-1}^{1-x}\rangle,\ 
\langle I_{i-1}^{1-x}\rangle 
\cdot\tilde v\cdot
b_{i-1}^\mu\langle I_{i-1}^x\rangle \bigr)\,\sigma
& \hbox{if $x_i=0$,} \\[5pt]
\bigl( \langle I_{i-1}^{1-x}\rangle 
\cdot\tilde u\cdot
b_{i-1}^\mu\langle I_{i-1}^x\rangle,\ 
\langle I_{i-1}^x\rangle\,b_{i-1}^{-\mu}
\cdot\tilde v\cdot 
\langle I_{i-1}^{1-x}\rangle\bigr)\,\sigma
& \hbox{if $x_i=1$.}
\end{array}$$
The next descendant $\tildew_{r+2-i}$ is now the concatenation of these entries and thus of the form
$$\begin{array}{ll}
\langle I_{i-1}^x\rangle\,b_{i-1}^{-\mu}
\cdot\tilde u\cdot 
\langle I_{i-1}^{1-x}\rangle
\cdot\tilde v\cdot
b_{i-1}^\mu\langle I_{i-1}^x\rangle
& \hbox{if $x_i=0$,} \\[5pt]
\langle I_{i-1}^{1-x}\rangle 
\cdot\tilde u\cdot
b_{i-1}^\mu
\langle I_{i-1}^x\rangle\,b_{i-1}^{-\mu}
\cdot\tilde v\cdot 
\langle I_{i-1}^{1-x}\rangle 
& \hbox{if $x_i=1$.}
\end{array}$$
In the case $x_i=0$ it follows that $\tildew_{r+2-i}$ is a product of  terms of the form $\langle I_{i-1}^x\rangle$ and/or $b_{i-1}^{-\mu}\langle I_{i-1}^{1-x}\rangle\,b_{i-1}^\mu$. In the case $x_i=1$ it is a product of  terms of the form $\langle I_{i-1}^{1-x}\rangle$ and/or $b_{i-1}^\mu\langle I_{i-1}^x\rangle\,b_{i-1}^{-\mu}$. In either case the lemma holds for $i-1$.

This finishes the proof of the induction step, and so the lemma follows for all $1\le i\le r$ by descending induction on~$i$.
\end{Proof}

\begin{Lem}\label{0MoreLem3}
The word $\tildew_r$ has the form $\langle I\rangle$ or $b_1^{-1} \langle I\rangle\,b_1$ or $b_1 \langle I\rangle\,b_1^{-1}$ for $I:= \{2,\ldots,r\}$.
\end{Lem}

\begin{Proof}
By Lemma \ref{1CombLem4} there exists $x\in\{0,1\}$ such that $I_1^x=\emptyset$. By Lemma \ref{1CombLem1} we then have $I_1^{1-x}=\{2,\ldots,r\}=I$. 
If $x=0$, Lemma \ref{0MoreLem2} implies that $\tildew_r$ has the form $\overline{b_1^{-1} \langle I\rangle\,b_1}$ or the form $\langle I\rangle\ \overline{b_1b_1^{-1} \langle I\rangle}$. But since $\tildew_r$ is minimal, it does not contain the subword $b_1b_1^{-1}$. Therefore $\tildew_r$ has the form $b_1^{-1} \langle I\rangle\,b_1$ or $\langle I\rangle$, as desired.
If $x=1$, Lemma \ref{0MoreLem2} implies that $\tildew_r$ has the form $\langle I\rangle\ \overline{b_1^{-1}b_1 \langle I\rangle}$ or the form $\overline{b_1 \langle I\rangle\,b_1^{-1}}$. Again by minimality $\tildew_r$ does not contain the subword $b_1^{-1}b_1$. It therefore has the form $\langle I\rangle$ or $b_1 \langle I\rangle\,b_1^{-1}$, as desired.
\end{Proof}

\begin{Lem}\label{0MoreLem4}
The word $\tildew_r$ has the form $\langle S^0\rangle$ or $b_1 \langle S^0\rangle\,b_1^{-1}$ or $\langle S^1\rangle$ or $b_1^{-1} \langle S^1\rangle\,b_1$.
\end{Lem}

\begin{Proof}
By Lemma \ref{0JIter} we have $J_{w_r} = \pi^r(J_{w_0}) = J_{w_0}$ and hence $1\not\in J_{w_r}$. Since $\tildew_r$ is a minimal word for the \OLGF-irreducible element~$w_r$, it therefore satisfies the conditions of Lemma \ref{0MoreLem1}.
Moreover, Lemma \ref{0MoreLem3} means that $\tildew_r$ contains the letters $b_1^{\pm1}$ at most in the first and last positions.
With Lemma \ref{0MoreLem1} it follows that $\tildew_r$ has the indicated form.
\end{Proof}

\begin{Lem}\label{0MoreLem5}
The word $\tildew_0$ has the form $\langle \{1\}\cup S^0\rangle$ or the form $\langle \{1\}\cup S^1\rangle$.
\end{Lem}

\begin{Proof}
For any $n\ge0$, each letter $b_i^{\pm1}$ of the word $\tildew_n$ bequeathes a letter \smash{$b_{\pi(i)}^{\pm1}$} to the word $\tildew_{n+1}$. By induction it follows that each letter $b_i^{\pm1}$ of the word $\tildew_0$ bequeathes a letter \smash{$b_{\pi^r(i)}^{\pm1}$}$\strut=b_i^{\pm1}$ to the word~$\tildew_r$. Thus $\tildew_0$ and $\tildew_r$ consist of the same letters, possibly rearranged. By Lemma \ref{0MoreLem4} these letters $b_i^{\pm1}$ either all satisfy $i\in \{1\}\cup S^0$ or all satisfy $i\in \{1\}\cup S^1$. Thus $\tildew_r$ and $\tildew_0$ have the indicated form.
\end{Proof}

%
%
%

\begin{Lem}\label{0MoreLem7}
For any \OLGF-irreducible element $w\in\GammaUx$ with $1\not\in J_w$,
the values $x_i$ for all $i\in J_w$ are equal.
\end{Lem}

\begin{Proof}
Apply the above constructions to $w_0 := w$. Then by Lemma \ref{0MoreLem5} the word $\tildew_0$ has the form $\langle \{1\}\cup S^x\rangle$ for some $x\in\{0,1\}$. By Lemma \ref{0LemIJ} we therefore have $J_w\subset \{1\}\cup S^x$ and hence $J_w\subset S^x$, as desired.
\end{Proof}

\begin{Prop}\label{0JCondLem}
For any \OLGF-irreducible $w\in\GammaUx$ the subset $J_w$ satisfies Condition \ref{0JCond}.
\end{Prop}

\begin{Proof}
Apply the above constructions to $w_0 := w$. Then by Lemma \ref{0JIter} for any $n\ge0$ we have $J_{w_n} = \pi^n(J_w)$. Thus if $1\not\in\pi^n(J_w)$, applying Lemma \ref{0MoreLem7} to $w_n$ implies that the values $x_i$ for all $i\in \pi^n(J_w)$ are equal. Therefore $J_w$ satisfies Condition \ref{0JCond}.
\end{Proof}


\subsection{Some rational functions and their denominators}
\label{sec0Rat}

To any subset $J\subset\{1,\ldots,r\}$ we now associate the following rational functions. Write the distinct elements of $J$ in ascending order $i_1<\ldots<i_k$ and set
\begin{eqnarray}
\UseTheoremCounterForNextEquation
\label{0PsiJDef}
\Psi_J &\!\!:=\!\!& \frac{1}{1-X^kY^r}\cdot \sum_{1\le j\le k} X^{j-1}Y^{i_j}
\qquad\hbox{and} \\[3pt]
\UseTheoremCounterForNextEquation
\label{0PhiJDef}
\Phi_J &\!\!:=\!\!& \frac{1}{1-2Y} \;+\; \frac{X-2}{\;1-2Y\;} \cdot \Psi_J.
\end{eqnarray}
For any $1\le k\le r$ we define
\UseTheoremCounterForNextEquation
\begin{equation}\label{0DkDef}
D_k\ :=\ \hbox{lowest common denominator of the $\Psi_J$ for all $J$ with $k = |J|$.}
\end{equation}

\begin{Prop}\label{0LCDAll}
For any $1\le k\le r$ we have
$$D_k \ =\ 
\left\{\begin{array}{ll}
(1-X^k Y^r) & \hbox{if $k<r$,} \\[5pt]
(1-XY) & \hbox{if $k=r$.} \\
\end{array}\right.$$
\end{Prop}

\begin{Proof}
By construction $D_k$ divides $1-X^k Y^r$. Conversely, for $J := \{1,\ldots,k\}$ we have
$$(1-X^k Y^r)\cdot \Psi_J \ = \sum_{1\le j\le k} X^{j-1}Y^j
\ =\ Y\cdot\frac{\ 1-X^kY^k}{1-XY}.$$
In the case $k<r$ the polynomials $1-X^kY^r$ and $Y(1-X^kY^k)$ are coprime, and hence $D_k=1-X^k Y^r$, as desired. In the case $k=r$ the only possible subset is $J=\{1,\ldots,r\}$, and then the above calculation shows that $\Psi_J = \frac{Y}{1-XY}$.
\end{Proof}


\subsection{Denominators of orbit length generating functions}
\label{sec0Denom}

\begin{Prop}\label{0IrredDenom}
For any \OLGF-irreducible element $w\in\GammaUx$ we have $\Phi_w = \Phi_{J_w}$.
\end{Prop}

\begin{Proof}
Let $i_1,\ldots,i_k$ be the distinct elements of~$J_w$, in any order. 
For all $n\ge0$ set $w_n=a_{\pi^n(i_1)}\cdots a_{\pi^n(i_k)}$. 
Then $w$ is $W$-conjugate to $w_0$ by Proposition \ref{0IrredConj}.
Also, by Lemma \ref{0Lem19} each $w_n$ is $W$-conjugate to $(w_{n+1},1)\,\sigma^{\mu_n}$, where $\mu_n=1$ if $1\in\pi^n(J_w)$ and $\mu_n=0$ otherwise. 
Since orbit length generating functions are invariant under conjugation, we have $\Psi_w=\Psi_{w_0}$, and with Proposition \ref{OLGFRecRels22} we deduce that
$$\Psi_{w_n} \ =\ 
\biggl\{\!\begin{array}{ll}
Y\Psi_1+Y\Psi_{w_{n+1}} & \hbox{if $1\not\in\pi^n(J_w)$,}\\[3pt]
Y+XY\Psi_{w_{n+1}} & \hbox{if $1\in\pi^n(J_w)$.}
\end{array}$$
As $\Psi_1=0$ by (\ref{F12}), the first of these formulas simplifies to $\Psi_{w_n} = Y\Psi_{w_{n+1}}$. By induction on $n$ it follows that 
$$\Psi_{w_0} \ =\ 
\! \sum_{0\le m<n \atop 1\in\pi^m(J_w)} \!\! X^{k_m}Y^{m+1} \;+\; X^{k_n}Y^n\Psi_{w_n},$$
where $k_m$ denotes the number of integers $0\le \ell<m$ such that $1\in\pi^\ell(J_w)$. Taking the limit in $\BZ[[X,Y]]$ we obtain
$$\Psi_w\ =\ \Psi_{w_0} \ =\ 
\! \sum_{m\ge0 \atop 1\in\pi^m(J_w)} \!\! X^{k_m}Y^{m+1}.$$
Since $\pi$ permutes the letters $1,\ldots,r$ transitively, and $J_w$ has cardinality~$k$, we have $k_{m+r}=k_m+k$ for all $m\ge0$. The last equality therefore implies that
$$\Psi_w \ =\ \frac{1}{1-X^kY^r} \cdot 
\! \sum_{0\le m<r \atop 1\in\pi^m(J_w)} \!\! X^{k_m}Y^{m+1}.$$
Moreover, for all $0\le m<r$ we have $1\in\pi^m(J_w)$ if and only if $m+1=\pi^{-m}(1)\in J_w$, and so
$$\Psi_w \ =\ \frac{1}{1-X^kY^r} \cdot \sum_{i\in J_w} X^{k_{i-1}}Y^i.$$
Also for all $i\in J_w$ we have $k_{i-1} = \bigl|\{m\in J_w\mid m<i\}\bigr|$.
Finally, since the last formula is independent of the order of $i_1,\ldots,i_k$, we may without loss of generality assume that $i_1<\ldots<\nobreak i_k$. Then for all $1\le j\le k$ we have
$k_{i_j-1} = \bigl|\{m\in J_w\mid m<i_j\}\bigr| = j-1$. Therefore $\Psi_w=\Psi_{J_w}$. By Proposition \ref{OLGFRecRels21} this implies that $\Phi_w = \Phi_{J_w}$, as desired.
\end{Proof}


\medskip
For any $1\le k\le r$ we now define
\UseTheoremCounterForNextEquation
\begin{equation}\label{0DkUxDef}
D_{\Ux,k}\ := \left[\;
\parbox{245pt}{lowest common denominator of the $\Psi_J$ for all $J$ with $k = |J|$ satisfying Condition \ref{0JCond}.}\;\right]
\end{equation}
By construction this is a divisor of the polynomial $D_k$ from (\ref{0DkDef}) and Proposition \ref{0LCDAll}.

\begin{Thm}\label{0MainThm}
The power series $\Phi_w \in 1+Y\BZ[[X,Y]]$ for all $w\in\GammaUx$ are rational functions with the lowest common denominator 
$$D_\Ux\ :=\ (1-2Y)\;\cdot\! \prod_{1\le k\le r}\!D_{\Ux,k}
\ \in\ 1+Y\BZ[X,Y].$$
\end{Thm}

\begin{Proof}
By Proposition \ref{0OLGFFin} and Theorem \ref{OLGFLinComb} the $\Phi_w$ for all $w\in\GammaUx$ are $\BZ[X,Y]$-linear combinations of the $\Phi_w$ for all \OLGF-irreducible elements. By Propositions \ref{0IrredExists} and \ref{0JCondLem} and \ref{0IrredDenom} the latter are precisely the $\Phi_J$ for all subsets $J\subset\{1,\ldots,r\}$ which satisfy Condition \ref{0JCond}. They are therefore rational functions, and in view of (\ref{0PhiJDef}) and (\ref{0DkUxDef}) their lowest common denominator is the least common multiple of the polynomials $(1-2Y)D_{\Ux,k}$ for all~$k$. 

The definition (\ref{0PsiJDef}) of $\Psi_J$ implies that each $D_{\Ux,k}$ divides $1-X^kY^r$. Thus $D_{\Ux,k}$ can be chosen congruent to $1\bmod Y$, and then $D_\Ux$ has the same property. 
Moreover, the polynomials $1-2Y$ and $1-X^kY^r$ for all $1\le k\le r$ are pairwise coprime, for instance because, viewed as polynomials in~$Y$, their zeros in an algebraic closure of $\BQ(X)$ are mutually distinct. Thus the least common multiple of all $(1-2Y)D_{\Ux,k}$ is~$D_\Ux$, and we are done.
\end{Proof}

\begin{Prop}\label{0Lem9}
\begin{enumerate}
\item[(a)] For all $\Ux$ we have $D_{\Ux,1}=D_1=1-XY^r$ and $D_{\Ux,r}=D_r=1-XY$.
\item[(b)] For $\Ux=(0,\ldots,0)$ or $(1,\ldots,1)$ we have $D_{\Ux,k}=D_k$ for all~$k$.
\end{enumerate}
\end{Prop}

\begin{Proof}
Assertion (a) follows from the fact that any subset $J$ of cardinality $1$ or $r$ satisfies Condition \ref{0JCond}. Assertion (b) follows from the fact that for these~$\Ux$, Condition \ref{0JCond} is satisfied for all~$J$.
\end{Proof}


In principle, the determination of the lowest common denominator $D_\Ux$ in Theorem \ref{0MainThm} is a finite combinatorial problem concerning the tuple~$\Ux$. The author does not (yet) know a simple direct description in general. However, we determined $D_\Ux$ in small cases using the computer algebra system Maple: see \cite{Pink2014Maple}. The outcome was that whenever $r\le 10$ and $x_2,\ldots,x_r$ are not all equal, then $\prod_{2\le k<r}D_{\Ux,k}=\nobreak1$ except in the following cases:
$$\begin{array}{|c|l|l|}
\hline
{\Large\strut} r & \, \prod_{2\le k<r}D_{\Ux,k} & \hbox{Conditions on $\Ux=(x_2,\ldots,x_r)$} \\[5pt]
\hline\hline
{\large\strut} 4 & (1-XY^2) & x_2=x_4\not=x_3 \\
\hline
{\large\strut} 6 & (1-X^2Y^6)(1-XY^2) & x_2=x_4=x_6\not=x_3=x_5 \\[-2pt]
{\large\strut} 6 & (1-XY^3)(1-X^2Y^3) & x_2=x_3=x_5=x_6\not=x_4 \\[-2pt]
{\large\strut} 6 & (1-XY^3) & x_2=x_5\not=x_3=x_6 \\[-2pt]
{\large\strut} 6 & (1-XY^2) & x_2=x_4=x_6 \longland x_3\not=x_5 \\
\hline
{\large\strut} 8 & (1-X^2Y^8)(1-X^3Y^8)(1-XY^2) & x_2=x_4=x_6=x_8\not=x_3=x_5=x_7 \\[-2pt]
{\large\strut} 8 & (1-XY^4)(1-X^2Y^4)(1-X^3Y^4) & x_2=x_3=x_4=x_6=x_7=x_8\not=x_5 \\[-2pt]
{\large\strut} 8 & (1-XY^4)(1-XY^2) & x_2=x_4=x_5=x_6=x_8\not=x_3=x_7 \\[-2pt]
{\large\strut} 8 & (1-XY^4) & x_2=x_6\not=x_4=x_8 \longland x_3=x_7 \\[-2pt]
{\large\strut} 8 & (1-XY^2) & x_2=x_4=x_6=x_8 \longland x_3\not=x_7 \\
\hline
{\large\strut} 9 & (1-X^2Y^9)(1-XY^3)(1-X^2Y^3) & x_2=x_3=x_5=x_6=x_8=x_9\not=x_4=x_7 \\[-2pt]
{\large\strut} 9 & (1-XY^3)(1-X^2Y^3) & x_2=x_3=x_5=x_6=x_8=x_9 \longland x_4\not=x_7 \\[-2pt]
{\large\strut} 9 & (1-X^2Y^9)(1-XY^3) & x_2=x_5=x_8\not=x_3=x_6=x_9 \longland x_4=x_7 \\[-2pt]
{\large\strut} 9 & (1-XY^3) & x_2=x_5=x_8\not=x_3=x_6=x_9 \longland x_4\not=x_7 \\
\hline
{\large\strut} 10 & \multicolumn{2}{l|}{(1-X^2Y^{10})(1-X^3Y^{10})(1-X^4Y^{10})(1-XY^2)} \\[-4pt]
{\large\strut} & \multicolumn{2}{r|}{x_2=x_4=x_6=x_8=x_{10}\not=x_3=x_5=x_7=x_9} \\
\hline
{\large\strut} 10 & \multicolumn{2}{l|}{(1-XY^5)(1-X^2Y^5)(1-X^3Y^5)(1-X^4Y^5)} \\[-4pt]
{\large\strut} & \multicolumn{2}{r|}{x_2=x_3=x_4=x_5=x_7=x_8=x_9=x_{10}\not=x_6} \\
\hline
{\large\strut} 10 & (1-XY^5) & x_2=x_7 \longland x_3=x_8 \longland x_4=x_9 \longland x_5=x_{10} \ \ \; \\[-2pt]
{\large\strut} && \hbox{but $x_2$, $x_3$, $x_4$, $x_5$ not all equal} \\
\hline
{\large\strut} 10 & (1-XY^2) & x_2=x_4=x_6=x_8=x_{10} \\[-2pt]
{\large\strut} && \hbox{but $x_3$, $x_5$, $x_7$, $x_9$ not all equal} \\
\hline
\end{array}$$



%
%
%
%
%
%
%

 \newpage


\hyphenation{polynomials}
\section{Iterated monodromy groups of quadratic polynomials: Pre-periodic case}
\label{S>0}
\hyphenation{poly-no-mi-als}


\subsection{The iterated monodromy group}
\label{sec1IMG}

Throughout this section we fix integers $r>s>0$ and a tuple $\Ux=(x_2,\ldots,x_r)$ of elements of $\{0,1\}$. Consider the elements $b_1,\ldots,b_r\in W$ defined by the recursion relations
\UseTheoremCounterForNextEquation
\begin{equation}\label{1bRecRels}
\left\{\begin{array}{lll}
b_1 &\!\!\!=\, \sigma, & \\[3pt]
b_{s+1} &\!\!\!=\, (b_r,b_s) & \hbox{if $x_{s+1}=0$,} \\[3pt]
b_{s+1} &\!\!\!=\, (b_s,b_r) & \hbox{if $x_{s+1}=1$,} \\[3pt]
b_i &\!\!\!=\, (b_{i-1},1) & \hbox{for all $i\not=1,s+1$ with $x_i=0$,} \\[3pt]
b_i &\!\!\!=\, (1,b_{i-1}) & \hbox{for all $i\not=1,s+1$ with $x_i=1$,} \\
\end{array}\right.
\end{equation}
and let $\GammaUx\subset W$ be the subgroup generated by them. Up to a change in notation, these are the generators and the subgroup studied by Bartholdi and Nekrashevych in \cite[\S4]{Bartholdi-Nekrashevych-2008}. Thus by \cite[Thm.\;5.1]{Bartholdi-Nekrashevych-2008} we have:

\begin{Thm}\label{1BarthNekrThm}
Let $f$ be any quadratic polynomial over $\BC$ and ${\eta\in\BC}$ be its unique critical point. Assume that $\eta,f(\eta),\ldots,f^r(\eta)$ are all distinct and that $f^{r+1}(\eta)=f^{s+1}(\eta)$. Then the iterated monodromy group of $f$ is $W$-conjugate to~$\GammaUx$ for a certain choice of the~$x_i$.
\end{Thm}

In the special case where $x_{s+1}=1$ and all other $x_i=0$ the above generators coincide with those studied in \cite[\S3]{Pink2013b}, but we do not care about them here. Instead, we will pay attention to the case where 
$\Ux=(0,\ldots,0)$,
that is, to the elements $a_1,\ldots,a_r$ defined by 
\UseTheoremCounterForNextEquation
\begin{equation}\label{1aRecRels}
\left\{\!\begin{array}{lll}
a_1      &\!\!\!=\, \sigma, & \\[3pt]
a_{s+1}  &\!\!\!=\, (a_r,a_s), & \\[3pt]
a_i      &\!\!\!=\, (a_{i-1},1) & \hbox{for all $i\not=1,s+1$.}\\
\end{array}\!\right.
\end{equation}
Also observe:

\begin{Prop}\label{1X1-X}
The group $\Gamma_{(x_2,\ldots,x_r)}$ is conjugate to the group $\Gamma_{(1-x_2,\ldots,1-x_r)}$ under~$W$.
\end{Prop}

\begin{Proof}
Same as that of Proposition \ref{0X1-X}, again with $w=(w,w)\,\sigma$.
\end{Proof}

\bigskip
The aim of this section is to show that the orbit length generating functions of all elements of $\GammaUx$ are rational and possess an explicit common denominator. 


\subsection{Finiteness}
\label{sec1Fin}

We begin with some preparations. 
Let $\pi$ denote the permutation of the set $\{1,\ldots,r\}$ defined by
\vskip-20pt
\UseTheoremCounterForNextEquation
\begin{equation}\label{1PiDef}
\pi(i)\ :=\ 
\left\{\begin{array}{ll}
s & \hbox{if $i=1$,}\\[3pt]
r & \hbox{if $i=s+1$,}\\[3pt]
i-1 & \hbox{otherwise.}
\end{array}\right.
\end{equation}
This induces a cyclic permutation of $\{1,\ldots,s\}$ and a cyclic permutation of $\{s+1,\ldots,r\}$. The recursion relations (\ref{1bRecRels}) express each $b_i$ in terms of $b_{\pi(i)}$, with $b_s$ thrown in for $i=s+1$ and taken out for $i=1$.

\medskip
Let $\Delta$ denote the subgroup of $\GammaUx$ that is generated by $b_1,\ldots,b_s$. The recursive description of these elements implies that $\Delta$ acts on the vertices of $T$ by changing only the last $s$ letters of a word, leaving the rest unchanged. Thus $\Delta$ is finite and acts faithfully on the subtree~$T_s$. In fact, one can easily show by induction that $\Delta$ maps isomorphically to the automorphism group of~$T_s$ and is therefore independent of~$\Ux$ (although the individual generators $b_1,\ldots,b_s$ depend on it). This characterization of $\Delta$ also implies:

\begin{Lem}\label{1Lem-1}
The only \OLGF-irreducible element of $\Delta$ is the identity element.
\end{Lem}

By contrast, repeated application of the recursion relations to $b_i$ for any $s<i\le r$ 
eventually leads back to~$b_i$. The two types of generators therefore play different roles in the arguments below. For instance, the letters $b_1,\ldots,b_s$ are not counted in the definition of the length below.

\begin{Lem}\label{1Lem0}
Every generator $b_i$ has order~$2$.
\end{Lem}

\begin{Proof}
Same as that of \cite[Prop.\;3.1.4]{Pink2013b}.
\end{Proof}

\begin{Def}\label{1LengthDef}
The \emph{length} $|w|$ of an element $w\in\GammaUx$ is the minimal number of letters from $\{b_{s+1},\ldots,b_r\}$ in a word over the alphabet $\{b_1,\ldots,b_r\}$ that represents~$w$. Any word representing $w$ 
with the minimal number of letters from $\{b_{s+1},\ldots,b_r\}$ is called a \emph{minimal word for~$w$}.
\end{Def}

Thus the elements of length $0$ of $\GammaUx$ are precisely those in~$\Delta$.

\begin{Lem}\label{1Lem1}
For any element $w=(u,v)\,\sigma^\mu \in\GammaUx$ we have $u,v\in\GammaUx$ and 
$${|uv| \le |u|+|v| \le |w|}.$$
\end{Lem}

\begin{Proof}
By the recursion relations (\ref{1bRecRels}), any letter $b_i$ for $s<i\le r$ in a minimal word for $w$ contributes precisely one letter $b_{\pi(i)}$ to a word representing precisely one of $u$,~$v$, and sometimes a letter $b_s$ which does not count towards the length. This implies the second inequality, and the first one follows directly from the definition of length.
\end{Proof}

\begin{Lem}\label{1Lem2}
For all $w\in\GammaUx$ and all $w'\in\Desc(w)$ we have $w'\in\GammaUx$ with $|w'|\le|w|$.
\end{Lem}

\begin{Proof}
By Definition \ref{DescDef} and iteration this follows from Lemma \ref{1LengthDef}.
\end{Proof}

\begin{Prop}\label{1OLGFFin}
Every element of $\GammaUx$ is \OLGF-finite.
\end{Prop}

\begin{Proof}
Any element of $\GammaUx$ of length $\ell$ can be written in the form $\delta_0 b_{i_1}\delta_1 \cdots b_{i_\ell}\delta_\ell$ with $\ell$ indices $s<i_j\le r$ and elements $\delta_j\in\Delta$. Since $\Delta$ is finite, it follows that $\GammaUx$ contains only finitely many elements of any given length. With Lemma \ref{1Lem2} this implies that $\Desc(w)$ is finite for any $w\in\GammaUx$, as desired.
\end{Proof}

\medskip
Combining Proposition \ref{1OLGFFin} with Theorem \ref{OLGFRat} we find that the orbit length generating functions of all elements of $\GammaUx$ are rational. By Theorem \ref{OLGFLinComb} the study of their denominators reduces to the case of \OLGF-irreducible elements. 
This case requires more preparations.


\subsection{Types and signs}
\label{sec1TypeSign}

\begin{Def}\label{1TypeDef}
An element $w\in\GammaUx$ is called 
\emph{of type} $I\subset\{1,\ldots,r\}$ if there exists a minimal word for $w$ which consists only of letters $b_i$ for $i\in I$.
\end{Def}

Note that this concerns a minimal word for~$w$, though in principle a minimal word might require a letter which some non-minimal word can do without. But this is intentional, because we use the notion of type as a secondary measure of complexity after the length.

\medskip
To any element $w\in\GammaUx$ we also associate the subset
\UseTheoremCounterForNextEquation
\begin{equation}\label{1JwDef}
J_w \ :=\ \{1\le i\le r\mid\sgn_i(w)=-1\}.
\end{equation}
To determine its relation with types we first observe:

\begin{Lem}\label{1LemSigns}
For all $1\le i\le r$ and $n\ge1$ we have
$$\sgn_n(b_i) \ =\ 
\left\{\!\begin{array}{rl}
-1 & \hbox{if $n=i\le s$,}\\[3pt]
-1 & \hbox{if $n\ge i>s$ and $n\equiv i\bmod\,(r-s)$,}\\[3pt]
 1 & \hbox{otherwise.}
\end{array}\right.$$
Thus for any fixed $w\in\GammaUx$, the value $\sgn_n(w)$ for $n>s$ depends only on $n\bmod(r-s)$.
\end{Lem}

\begin{Proof} 
Same as that of \cite[Prop.\;3.1.1]{Pink2013b}.
\end{Proof}

\begin{Lem}\label{1LemIJ}
If $w\in\GammaUx$ is of type~$I$, then $J_w\subset I$.
\end{Lem}

\begin{Proof}
Lemma \ref{1LemSigns} implies that $\sgn_i(b_i)=-1$ and ${\sgn_i(b_j)=1}$ whenever $i\not=j$.
\end{Proof}


\subsection{Properties of \OLGF-irreducible elements}
\label{sec1Prop}

\begin{Lem}\label{1Lem3}
Any \OLGF-irreducible element $w\in\GammaUx$ has a unique first descendant $w'$ which is \OLGF-irreducible with $|w'|=|w|$. Moreover $w$ is either $W$-conjugate to $(w',1)\,\sigma$, or equal to $(w',\delta)$ or $(\delta,w')$ with $\delta\in\Delta$ where $w'$ and $\delta$ are the first descendants of~$w$.
\end{Lem}

\begin{Proof}
Suppose first that $w=(u,v)\,\sigma$. Then $w$ is $W$-conjugate to $(uv,1)\,\sigma$, and $uv$ is the unique first descendant of~$w$. Thus the assumption $w\in\Desc(w)$ means that $w$ is equal to or a descendant of~$uv$. On the one hand this implies that $uv$ is a descendant of itself; hence $uv$ is \OLGF-irreducible. On the other hand it implies by Lemma \ref{1Lem2} that $|w|\le|uv|\le|w|$ and hence $|uv|=|w|$, and we are done with $w':=uv$.

Suppose now that $w=(u,v)$, so that $u$ and $v$ are the first descendants of~$w$. Then the assumption $w\in\Desc(w)$ means that $w$ is equal to, or a descendant of, one of $u$, $v$; let us call it~$w'$. On the one hand this implies that $w'$ is a descendant of itself; hence $w'$ is \OLGF-irreducible. On the other hand it implies by Lemma \ref{1Lem2} that $|w|\le|w'|\le|w|$ and hence $|w'|=|w|$. Plugging this into the inequality $|u|+|v|\le|w|$ from Lemma \ref{1Lem1}, we now deduce that the other entry of $(u,v)$ has length $0$ and therefore lies in~$\Delta$. Calling it $\delta$, we then have $w=(w',\delta)$ or $w=(\delta,w')$. Finally this makes $w'$ unique unless $|w|=0$. But in that case $w=1=(1,1)$ by Lemma \ref{1Lem-1} and hence $w'=\delta=1$ is again unique, and we are done.
\end{Proof}

\medskip
For the following arguments we fix a \OLGF-irreducible element $w_0\in\GammaUx$ and construct a sequence of \OLGF-irreducible elements $w_n\in\GammaUx$ by defining each $w_{n+1}$ as the first descendant of $w_n$ furnished by Lemma \ref{1Lem3}. 

\begin{Lem}\label{1LemSeqPer}
The sequence $w_0,w_1,\ldots$ is periodic.
\end{Lem}

\begin{Proof}
Repeated application of Lemma \ref{1Lem3} shows that the descendants of $w_0$ are precisely the elements $w_n$ for $n\ge1$ and perhaps some elements of~$\Delta$. Since by assumption $w_0$ is a descendant of itself, we must have $w_{n_0}=w_0$ for some $n_0\ge1$. Then the construction implies that $w_{n+n_0}=w_n$ for all ${n\ge0}$. 
\end{Proof}

\begin{Lem}\label{1LemType1}
Consider any $n\ge0$, and suppose that $w_n$ is of type~$I$. Then $w_{n+1}$ is of type $\pi(I)$. If moreover $s+1\not\in I$ or $1\not\in J_{w_n}$, 
then $w_{n+1}$ is of type $\pi(I)\setminus\{s\}$.
\end{Lem}

\begin{Proof}
Set $\ell:=|w_n|$ and write $w_n$ as a minimal word $w_n=b_{i_1}\cdots b_{i_k}$ with all $i_j\in I$. Then precisely $\ell$ of these letters lie in $\{b_{s+1},\ldots,b_r\}$. Write $w_n=(u,v)\sigma^\mu$ and use the recursion relations (\ref{1bRecRels}), but no other relations, to obtain words representing $u$ and~$v$. Then precisely $\ell$ of the letters of both words lie in $\{b_{s+1},\ldots,b_r\}$. Since $w_{n+1}$ is equal to $u$, $v$, or $uv$, and itself of length $\ell$ by Lemma \ref{1Lem3}, this results in a minimal word representing~$w_{n+1}$.

By construction, any letter $b_{i_j}\not=b_1,b_{s+1}$ contributes at most one letter $b_{i_j-1}$ to the word representing~$w_{n+1}$. Any letter $b_{i_j}=b_{s+1}$ contributes at most the letters $b_r$ and~$b_s$, and any letter $b_{i_j}=b_1$ contributes nothing. This shows that $w_{n+1}$ is of type
$$I'\ := \left\{\begin{array}{ll}
\{i-1\mid1,s+1\not=i\in I\} 
& \hbox{if $s+1\not\in I$,} \\[5pt]
\{i-1\mid1,s+1\not=i\in I\} \cup \{r,s\}
& \hbox{if $s+1\in I$.}
\end{array}\right.$$

If $s+1\not\in I$, the definition (\ref{1PiDef}) of $\pi$ implies that $I'=\pi(I)\setminus\{s\}$, and we are done. For the rest of the proof we therefore assume that $s+1\in I$.

If $1\in J_{w_n}$, we have $1\in I$ by Lemma \ref{1LemIJ}. Since also $s+1\in I$, the definition of $\pi$ now implies that $I'=\pi(I)$, and again we are done. For the rest of the proof we therefore assume that 
$1\not\in J_{w_n}$.

Then $\sgn_1(w_n)=+1$ and hence $w_n=(u,v)$. Also, one of its entries is~$w_{n+1}$, and the word representing it contains all $\ell$ occurrences of letters in $\{b_{s+1},\ldots,b_r\}$ that result from letters $b_{i_j}\in\{b_{s+1},\ldots,b_r\}$. In particular, the word representing $w_{n+1}$ contains all occurrences of the letter $b_r$ resulting from a letter $b_{i_j}=b_{s+1}$. Each such $b_{i_j}$ contributes a letter $b_s$ to the other entry of $(u,v)$. Since the letter $b_s$ does not arise in any other way from the recursion relations, it follows that the letter $b_s$ does not occur in the word representing~$w_{n+1}$. Therefore $w_{n+1}$ is of type $I'\setminus\{s\}$. But the definition (\ref{1PiDef}) of $\pi$ implies that $I'\setminus\{s\}=\pi(I)\setminus\{s\}$, so we are done.
\end{Proof}

\medskip
Next we fix a subset $I_0\subset\{1,\ldots,r\}$ of minimal cardinality such that $w_0$ is of type~$I_0$. For every $n\ge0$ we set $I_n := \pi^n(I_0)$. 

\begin{Lem}\label{1LemType2}
For every $n\ge0$, the set $I_n\subset\{1,\ldots,r\}$ is a subset of minimal cardinality such that $w_n$ is of type~$I_n$. 
\end{Lem}

\begin{Proof}
By induction on~$n$, Lemma \ref{1LemType1} implies that $w_n$ is of type $I_n$ for all $n\ge0$. Suppose that for some $n\ge0$ there exists a subset $I'_n \subset \{1,\ldots,r\}$ with $|I'_n|<|I_n|$ such that $w_n$ is of type~$I_n'$. Then again by Lemma \ref{1LemType1}, the element $w_{n'}$ is of type $\pi^{n'-n}(I'_n)$ for every $n'\ge n$. By Lemma \ref{1LemSeqPer} we can choose $n'\ge n$ such that $w_{n'}=w_0$. Then $w_0$ is of type $\pi^{n'-n}(I'_n)$ with $|\pi^{n'-n}(I'_n)| = |I'_n|<|I_n|=|I_0|$, contradicting the minimality of~$|I_0|$. 
\end{Proof}

\begin{Lem}\label{1LemType3}
For any $n\ge0$ we have $1\in J_{w_n}$ if and only if $1\in I_n$.
\end{Lem}

\begin{Proof}
The `only if' part follows from Lemmas \ref{1LemIJ} and \ref{1LemType2}. For  the `if' part suppose that $1\not\in J_{w_n}$. Then $w_{n+1}$ is of type $\pi(I_n)\setminus\{s\}$ by Lemma \ref{1LemType1}. The minimality of $I_{n+1}=\pi(I_n)$ from Lemma \ref{1LemType2} then implies that $\pi(I_n)\setminus\{s\} = \pi(I_n)$. Thus $s\not\in\pi(I_n)$, and hence $1=\pi^{-1}(s)\not\in I_n$, proving the converse.
\end{Proof}

\begin{Lem}\label{1LemType4}
For any $n\ge0$ with $1\not\in I_n$, the values $x_i$ are equal for all $i\in I_n$, and $w_n$ is equal to $(w_{n+1},b_s^{\nu_n})$ or $(b_s^{\nu_n},w_{n+1})$ for $\nu_n\in\BZ$ such that $\sgn_{s+1}(w_n)=(-1)^{\nu_n}$. 
\end{Lem}

\begin{Proof}
Assume that $1\not\in I_n$ and abbreviate $I_n^x := \{i\in I_n\mid x_i=x\}$ for all $x\in\{0,1\}$. 
Write $w_n=b_{i_1}\cdots b_{i_k}$ as a minimal word with all $i_j\in I_n$. Since $1\not\in I_n$, the recursion relations (\ref{1bRecRels}) show that all factors have the form 
$$\left\{\begin{array}{ll}
(b_r,b_s) & \hbox{if $s+1\in I_n^0$,} \\[3pt]
(b_s,b_r) & \hbox{if $s+1\in I_n^1$,} \\[3pt]
(b_{i-1},1) & \hbox{for $i\in I_n^0\setminus\{s+1\}$,} \\[3pt]
(1,b_{i-1}) & \hbox{for $i\in I_n^1\setminus\{s+1\}$.}
\end{array}\right.$$
By Lemma \ref{1Lem3} we have $w_n=(w_{n+1},\delta)$ or $(\delta,w_{n+1})$ for some $\delta\in\Delta$. Moreover, in the proof of Lemma \ref{1LemType1} we have seen that the expansions in the above list yield a minimal word for~$w_{n+1}$. Set $x:=0$ if $w_n=(w_{n+1},\delta)$, and $x:=1$ otherwise. Then the above list implies that the resulting word for $w_{n+1}$ is a product of certain $b_j$ for $j\in \pi(I_n^x) \cup \{s\}$. But in the proof of Lemma \ref{1LemType1} we have already seen that in this case all occurrences of $b_s$ must go into~$\delta$. Thus $w_{n+1}$ is of type $\pi(I_n^x)$.

Now the minimality in Lemma \ref{1LemType2} implies that the inclusion $\pi(I_n^x) \subset \pi(I_n)=I_{n+1}$ is an equality. Thus $I_n^x=I_n$, proving the first assertion. Plugging this back into the above list now shows that the only non-trivial factors going into $\delta$ are the $b_s$ arising from all $b_{i_j}=b_{s+1}$. Thus if $\nu_n$ denotes the number of factors $b_{i_j}=b_{s+1}$, we have $\delta=b_s^{\nu_n}$. But then Lemma \ref{1LemSigns} shows that $\sgn_{s+1}(w_n)=(-1)^{\nu_n}$, and we are done.
\end{Proof}

\begin{Lem}\label{1LemInduc}
For any $n\ge0$ the element $w_n$ is $W$-conjugate to
$$\begin{array}{ll}
(w_{n+1},1)\,\sigma & \hbox{if $1\in J_{w_n}$,} \\[3pt]
(w_{n+1},b_s) & \hbox{if $1\not\in J_{w_n}$ and $s+1\in J_{w_n}$,} \\[3pt]
(w_{n+1},1) & \hbox{if $1\not\in J_{w_n}$ and $s+1\not\in J_{w_n}$.}
\end{array}$$
\end{Lem}

\begin{Proof}
If $1\in J_{w_n}$, that is, if $\sgn_1(w_n)=-1$, this follows from Lemma \ref{1Lem3}. Otherwise we have $1\not\in I_n$ by Lemma \ref{1LemType3}, and so the remaining cases follow from Lemma \ref{1LemType4}.
\end{Proof}

\begin{Lem}\label{1LemFourSgn}
For every $n\ge0$ we have:
\begin{enumerate}
\item[(a)] For all $i\ge2$ with $i\not=s+1$ we have $\sgn_i(w_n)=\sgn_{i-1}(w_{n+1})$.
\item[(b)] For all $1\le i\le s$ we have $\sgn_i(w_n)=-1$ if and only if \ $i\in I_n$.
\item[(c)] For all $1\le i\le r$ we have $\sgn_i(w_n)=\sgn_{\pi(i)}(w_{n+1})$.
\item[(d)] If $\sgn_1(w_n)=-1$, then $\sgn_{s+1}(w_n)=-1$.
\end{enumerate}
\end{Lem}

\begin{Proof}
For all $i\ge2$ with $i\not=s+1$ we have $\sgn_i(\sigma)=1$ and $\sgn_{i-1}(b_s)=1$ by Lemma \ref{1LemSigns}. Thus in each of the cases in Lemma \ref{1LemInduc}, the conjugation invariance and the recursion relations for signs imply that $\sgn_i(w_n)=\sgn_{i-1}(w_{n+1})$, proving (a).

Next we prove (b) simultaneously for all $n$ by induction on~$i$. For $i=1$ the assertion already holds by Lemma \ref{1LemType3}. If $i>1$, by (a) we have $\sgn_i(w_n)=-1$ if and only if $\sgn_{i-1}(w_{n+1})=-1$. By the induction hypothesis this is equivalent to $i-1\in I_{n+1}$, in other words to $\pi(i)=i-1\in I_{n+1}=\pi(I_n)$, and hence to $i\in I_n$, finishing the induction step.

Assertion (c) for $i\not=1,s+1$ is the same as (a). For $i=1$ by (b) we have 
$\sgn_1(w_n)=-1$ if and only if $1\in I_n$ if and only if $s=\pi(1) \in \pi(I_n)=I_{n+1}$, which again by (b) is equivalent to $\sgn_{\pi(1)}(w_{n+1})=-1$. This proves (c) for $i=1$. For $i=s+1$ by the periodicity in Lemma \ref{1LemSigns} combined with (a) we have $\sgn_{s+1}(w_n) = \sgn_{r+1}(w_n) = \sgn_r(w_{n+1}) = \sgn_{\pi(s+1)}(w_{n+1})$. This proves (c) in all cases.

To show (d) we repeat the argument for (a) with $s+1$ in place of~$i$. Again we have $\sgn_{s+1}(\sigma)=1$, and since now $\sgn_1(w_n)=-1$, the first case of Lemma \ref{1LemInduc} implies that $\sgn_{s+1}(w_n)=\sgn_s(w_{n+1})$. By (c) this is equal to $\sgn_{\pi(1)}(w_{n+1}) = \sgn_1(w_n) = -1$, as desired. Thus everything is proved.
\end{Proof}

\begin{Lem}\label{1JIter}
For every $n\ge0$ we have $J_{w_n} = \pi^n(J_{w_0})$.
\end{Lem}

\begin{Proof}
This follows by induction from Lemma \ref{1LemFourSgn} (c).
\end{Proof}

\medskip
Now consider the following conditions on a subset $J\subset\{1,\ldots,r\}$:

\begin{Conds}\label{1JConds}
For all $n\ge0$,
\begin{enumerate}
\item[(a)] if $1\in\pi^n(J)$, then $s+1\in\pi^n(J)$.
\item[(b)] if $1\not\in\pi^n(J)$, then the values $x_i$ are equal for all $i\in\pi^n(J)$.
\end{enumerate}
\end{Conds}

\begin{Prop}\label{1JCondLem}
For any \OLGF-irreducible $w\in\GammaUx$ the subset $J_w$ satisfies Conditions \ref{1JConds}.
\end{Prop}

\begin{Proof}
Apply the above with $w_0:=w$, and consider any $n\ge0$. If $1\in\pi^n(J_{w_0})$, by Lemma \ref{1JIter} we have $\sgn_1(w_n)=-1$. By Lemma \ref{1LemFourSgn} (d) this implies that $\sgn_{s+1}(w_n)=-1$ and therefore $s+1\in\pi^n(J_{w_0})$, proving the condition \ref{1JConds} (a). 
By contrast, if $1\not\in\pi^n(J_{w_0})$, then $1\not\in I_n$ by Lemma \ref{1LemType3}. From Lemma \ref{1LemType4} it then follows that the values $x_i$ are equal for all $i\in I_n$. But Lemmas \ref{1JIter} and \ref{1LemIJ} together imply that $\pi^n(J_{w_0})\subset I_n$, so in particular the values $x_i$ are equal for all $i\in\pi^n(J_{w_0})$, proving the condition \ref{1JConds} (b). 
\end{Proof}

\begin{Lem}\label{1JCondsaEq}
Condition \ref{1JConds} (a) is equivalent to:
\begin{enumerate}
\item[(a\,$'$)] For all $i\in J$ with $i\le s$ and all $s<j\le r$ with $i\equiv j\bmod(s,r-s)$ we have $j\in J$.
\end{enumerate}
\end{Lem}

\begin{Proof}
For any $n\ge0$ we have $1\in\pi^n(J)$ if and only if $\pi^{-n}(1)\in J$, and $\pi^{-n}(1)$ is the unique integer $1\le i\le s$ with $i\equiv n+1\bmod(s)$. Similarly, we have $s+1\in\pi^n(J)$ if and only if $\pi^{-n}(s+1)\in J$, where $\pi^{-n}(s+1)$ is the unique integer $s<j\le r$ with $j\equiv s+n+1\bmod(r-s)$. Given $i$ and~$j$, the conditions on $n$ just stated are $n+s+1\equiv i\bmod(s)$ and $n+s+1\equiv j\bmod(r-s)$, so they are satisfied for some $n$ if and only if $i\equiv j\bmod(s,r-s)$. Now the equivalence follows.
\end{Proof}


\subsection{Conjugacy classes of \OLGF-irreducible elements}
\label{sec1Conj}

\begin{Lem}\label{1Lem19}
Consider any distinct indices $i_1,\ldots,i_k \in \{1,\ldots,r\}$, in any order, and set $J:=\{i_1,\ldots,i_k\}$. If\/ $1\in J$ assume that also $s+1\in J$. Then $a_{i_1}\cdots a_{i_k}$ is $W$-conjugate to 
$$\begin{array}{ll}
(a_{\pi(i_1)}\cdots a_{\pi(i_k)},1)\,\sigma 
& \hbox{if $1\in J$ and $s+1\in J$,} \\[3pt]
(a_{\pi(i_1)}\cdots a_{\pi(i_k)},a_s)
& \hbox{if $1\not\in J$ and $s+1\in J$,} \\[3pt]
(a_{\pi(i_1)}\cdots a_{\pi(i_k)},1)
& \hbox{if $1\not\in J$ and $s+1\not\in J$.}
\end{array}$$
\end{Lem}

\begin{Proof}
If $1\not\in J$ and $s+1\not\in J$, the recursion relations (\ref{1aRecRels}) imply that 
$$a_{i_1}\cdots a_{i_k} \ =\  
(a_{i_1-1},1)\cdots(a_{i_k-1},1) \ =\ 
(a_{\pi(i_1)}\cdots a_{\pi(i_k)},1),$$
and the assertion follows. If $1\not\in J$ and $s+1\in J$, let $j$ be the unique index with $i_j=s+1$. Then the recursion relations (\ref{1aRecRels}) imply that 
$$\begin{array}{rl}
a_{i_1}\cdots a_{i_k} \ =&
(a_{i_1-1},1)\cdots(a_{i_{j-1}-1},1)
\cdot(a_r,a_s) \cdot
(a_{i_{j+1}-1},1)\cdots(a_{i_k-1},1) \\[3pt]
\ =& (a_{\pi(i_1)}\cdots a_{\pi(i_k)},a_s),
\end{array}$$
and again the assertion follows. So assume that $1\in J$ and $s+1\in J$, and let $\ell$ and $j$ be the unique indices with $i_\ell=1$ and $i_j=s+1$. If $\ell<j$, the recursion relations (\ref{1aRecRels}) imply that 
$$\begin{array}{rl}
\kern-15pt a_{i_1}\cdots a_{i_k} \ =&
(a_{i_1-1},1)\cdots(a_{i_{\ell-1}-1},1)
\, \sigma \,
(a_{i_{\ell+1}-1},1)\cdots(a_{i_{j-1}-1},1)
\, (a_r,a_s) \,
(a_{i_{j+1}-1},1)\cdots(a_{i_k-1},1) \\[3pt]
\ =&
\bigl(a_{\pi(i_1)}\cdots a_{\pi(i_{\ell-1})}a_s\,,\,
a_{\pi(i_{\ell+1})}\cdots a_{\pi(i_{j-1})}
a_r
a_{\pi(i_{j+1})}\cdots a_{\pi(i_k)}\bigr) 
\,\sigma \\[3pt]
\ =&
(a_{\pi(i_1)}\cdots a_{\pi(i_{\ell})}\,,\,
a_{\pi(i_{\ell+1})}\cdots a_{\pi(i_k)}) 
\,\sigma
\end{array}$$
which is conjugate to $(a_{\pi(i_1)}\cdots a_{\pi(i_k)},1)\,\sigma$, as desired. If $\ell>j$, the same kind of calculation yields
$$\begin{array}{rl}
\kern-15pt a_{i_1}\cdots a_{i_k} \ =&
(a_{i_1-1},1)\cdots(a_{i_{j-1}-1},1)
\, (a_r,a_s) \,
(a_{i_{j+1}-1},1)\cdots(a_{i_{\ell-1}-1},1)
\, \sigma \,
(a_{i_{\ell+1}-1},1)\cdots(a_{i_k-1},1) \\[3pt]
\ =&
\bigl(a_{\pi(i_1)}\cdots a_{\pi(i_{j-1})}
a_r
a_{\pi(i_{j+1})}\cdots a_{\pi(i_{\ell-1})}\,,\,
a_s
a_{\pi(i_{\ell+1})}\cdots a_{\pi(i_k)}\bigr) 
\,\sigma \\[3pt]
\ =&
(a_{\pi(i_1)}\cdots a_{\pi(i_{\ell-1})}\,,\,
a_{\pi(i_{\ell})}\cdots a_{\pi(i_k)}) 
\,\sigma
\end{array}$$
which is again conjugate to $(a_{\pi(i_1)}\cdots a_{\pi(i_k)},1)\,\sigma$, as desired.
\end{Proof}

\begin{Lem}\label{1Lem20}
For each $1\le i\le r$ the element $b_i$ is conjugate to $a_i$ under~$W$.
\end{Lem}

\begin{Proof}
Let $a'_1,\ldots,a'_r$ denote the generators used in \cite[\S3]{Pink2013b}. Then the recursion relations (\ref{1bRecRels}) and the equivalence (a)$\Leftrightarrow$(b) of \cite[Thm.\;3.4.1]{Pink2013b} imply that each $b_i$ is individually conjugate to $a'_i$ under~$W$. Since the elements $a_1,\ldots,a_r$ are a special case of the elements $b_1,\ldots,b_r$,
it follows that each $b_i$ is individually conjugate to $a_i$ under~$W$. 
\end{Proof}

\begin{Prop}\label{1IrredConj}
Consider any \OLGF-irreducible element $w\in\GammaUx$. Let $i_1,\ldots,i_k$ be the distinct elements of $J_w$ in any order. Then $w$ is conjugate to $a_{i_1}\cdots a_{i_k}$ under~$W$.
\end{Prop}


\begin{Proof}
By \cite[Lemma\;1.3.3]{Pink2013b} it suffices to show that the restrictions $w|_{T_n}$ and $a_{i_1}\cdots a_{i_k}|_{T_n}$ are conjugate in the automorphism group of $T_n$ for every $n\ge0$. We will achieve this by induction on~$n$. For $n=0$ the assertion is trivially true, so assume that $n>0$ and that the assertion is universally true for the restrictions to $T_{n-1}$.

We apply the above constructions to $w_0 := w$. Then from Lemma \ref{1JIter} we know that $\pi(i_1),\ldots,\pi(i_k)$ are the distinct elements of $J_{w_1}$. By the induction hypothesis the restriction $w_1|_{T_{n-1}}$ is therefore conjugate to $a_{\pi(i_1)}\cdots a_{\pi(i_k)}|_{T_{n-1}}$ under the automorphism group of~$T_{n-1}$. 
Since $J_{w_0}$ satisfies Condition \ref{1JConds} (a), by Lemma \ref{1Lem19} it follows that $a_{i_1}\cdots a_{i_k}|_{T_n}$ is conjugate to 
$$\begin{array}{ll}
(w_1,1)\,\sigma|_{T_n} & \hbox{if $1\in J_{w_0}$,} \\[3pt]
(w_1,a_s)|_{T_n} & \hbox{if $1\not\in J_{w_0}$ and $s+1\in J_{w_0}$,} \\[3pt]
(w_1,1)|_{T_n} & \hbox{if $1\not\in J_{w_0}$ and $s+1\not\in J_{w_0}$,}
\end{array}$$
under the automorphism group of~$T_n$. Moreover, by Lemma \ref{1Lem20} the element $b_s$ is conjugate to $a_s$ under~$W$. Comparing the above cases with the respective cases in Lemma \ref{1Lem19} thus implies that $w_0|T_n$ is conjugate to $a_{i_1}\cdots a_{i_k}|_{T_n}$ under the automorphism group of~$T_n$. This finishes the induction step and thereby the proof.
\end{Proof}


\begin{Prop}\label{1IrredExists}
For any subset $J\subset\{1,\ldots,r\}$ satisfying Conditions \ref{1JConds} there exists a \OLGF-irreducible element $w\in\GammaUx$ with $J = J_w$.
\end{Prop}

\begin{Proof}
Consider any integer $n\ge0$. For the purpose of this proof we call any element of $\GammaUx$ of the form $b_{i_1}\cdots b_{i_k}$, where $i_1,\ldots,i_k$ are the distinct elements of $\pi^n(J)$ in any order, \emph{strictly of type $\pi^n(J)$}. (This is actually a minimal word; hence the element is of type $\pi^n(J)$ in the sense of Definition \ref{1TypeDef}, but we will not use that fact.) We claim that any element that is strictly of type $\pi^n(J)$ possesses a first descendant which is strictly of type $\pi^{n+1}(J)$.

Granting this, by induction on $n$ it follows that for any $n\ge1$, any element that is strictly of type $J$ possesses a descendant which is strictly of type $\pi^n(J)$. Since $\pi$ is a permutation of finite order, we deduce that any element that is strictly of type $J$ possesses a descendant which is again strictly of type~$J$. As there are only finitely many elements that are strictly of type $J$, and being a descendant is a transitive relation, it follows that some element $w$ that is strictly of type $J$ must be its own descendant. This element is therefore \OLGF-irreducible. Finally, writing $w=b_{i_1}\cdots b_{i_k}$ where $i_1,\ldots,i_k$ are the distinct elements of~$J$, Lemma \ref{1LemSigns} implies that $J = J_w$, as desired.

To prove the claim consider $w := b_{i_1}\cdots b_{i_k}$ where $i_1,\ldots,i_k$ are the distinct elements of $\pi^n(J)$. Suppose first that $1\not\in\pi^n(J)$. 
Then by Condition \ref{1JConds} (b) the values $x_i$ are equal for all $i\in\pi^n(J)$. If this common value is~$0$, the same calculations as in the proof of Lemma \ref{1Lem19} show that $w$ is equal to $(b_{\pi(i_1)}\cdots b_{\pi(i_k)},1)$ or $(b_{\pi(i_1)}\cdots b_{\pi(i_k)},b_s)$. If the common value is~$1$, then $w$ is equal to $(1,b_{\pi(i_1)}\cdots b_{\pi(i_k)})$ or $(b_s,b_{\pi(i_1)}\cdots b_{\pi(i_k)})$. In all these cases $w$ has the first descendant $b_{\pi(i_1)}\cdots b_{\pi(i_k)}$, which is strictly of type $\pi^{n+1}(J)$.

Suppose now that $1\in\pi^n(J)$. Then $\sgn_1(w)=-1$ by Lemma \ref{1LemSigns}; hence $w$ has the form $w=(u,v)\,\sigma$. By Condition \ref{1JConds} (a) we now also have $s+1\in\pi^n(J)$. By the recursion relations (\ref{1bRecRels}), any factor $b_{i_j}\not=b_1,b_{s+1}$ of $w=b_{i_1}\cdots b_{i_k}$ contributes precisely one factor $b_{i_j-1}= b_{\pi(i_j)}$ to the product~$uv$. The factor $b_{i_j}=b_{s+1}$ contributes the factors $b_r=b_{\pi(s+1)}$ and~$b_s=b_{\pi(1)}$, and the factor $b_{i_j}=b_1$ contributes nothing. Together this shows that $uv$ is a product of the elements 
$b_{\pi(i_1)},\ldots,b_{\pi(i_k)}$ in some order. It is therefore strictly of type $\pi^{n+1}(J)$, as desired.
\end{Proof}

\begin{Prop}\label{1IrredConjToSplit}
Any \OLGF-irreducible element $w$ of $\GammaUx$ is $W$-conjugate to a \OLGF-irreducible element of $\Gamma_{(0,\ldots,0)}$.
\end{Prop}

\begin{Proof}
By Proposition \ref{1IrredConj} it is conjugate to $a_{i_1}\cdots a_{i_k} \in \Gamma_{(0,\ldots,0)}$, where $i_1,\ldots,i_k$ are the distinct elements of~$J_w$ in any order. But the same argument as in the proof of Proposition \ref{1IrredExists} shows that for some order, the element $a_{i_1}\cdots a_{i_k}$ is \OLGF-irreducible.
\end{Proof}


\subsection{Some rational functions and their denominators}
\label{sec1Rat}

To any subset $J\subset\{1,\ldots,r\}$ we associate the following power series. For all $m\ge0$ let $\ell_m$ denote the number of integers $0\le i<m$ such that $1\in\pi^i(J)$. Set
\begin{eqnarray}
\UseTheoremCounterForNextEquation
\label{1PsiJDef}
\Psi_J &\!\!:=\!\!& \sum_{m\ge0 \atop 1\in\pi^m(J)} \!\! X^{\ell_m}Y^{m+1}
\; + \kern-10pt \sum_{m\ge0 \atop 1\not\in\pi^m(J)\ni s+1} 
\kern-13pt X^{\ell_m}Y^{m+s+1},
\rlap{\qquad\hbox{and}} \\[3pt]
\UseTheoremCounterForNextEquation
\label{1PhiJDef}
\Phi_J &\!\!:=\!\!& \frac{1}{1-2Y} \;+\; \frac{X-2}{\;1-2Y\;} \cdot \Psi_J.
\end{eqnarray}
Abbreviate $p:= {\gcd(s,r-s)}$ and $q := \frac{r-s}{p}$.

\begin{Lem}\label{1PsiRat}
Both $\Psi_J$ and $\Phi_J$ are rational functions. More precisely, 
with $\ell := \bigl|\{{i\in J}\mid {i\le s}\}\bigr|$ we have
$$(1-X^{\ell q}Y^{sq})\cdot \Psi_J \ \in\ \BZ[X,Y].$$
\end{Lem}

\begin{Proof}
Since $\pi$ permutes the letters $1,\ldots,s$ transitively, the condition $1\in\pi^m(J)$ depends only on $m\bmod s$, and we have $\ell_{m+s}=\ell_m+\ell$ for all $m\ge0$. Also, since $\pi$ permutes the letters $s+1,\ldots,r$ transitively, the condition $s+1\in\pi^m(J)$ depends only on $m\bmod(r-s)$. The definition of $p$ and $q$ implies that $sq=\mathop{\rm lcm}(s,r-s)$. Thus in both sums in (\ref{1PsiJDef}), the terms with $m+sq$ in place of $m$ are obtained on multiplying the terms for $m$ by $X^{\ell q}Y^{sq}$. Therefore $(1-X^{\ell q}Y^{sq})\cdot \Psi_J \in \BZ[X,Y]$, and $\Psi_J$ is rational. By (\ref{1PhiJDef}) so is~$\Phi_J$.
\end{Proof}

\medskip
For any $0\le\ell\le s$ we define:
\UseTheoremCounterForNextEquation
\begin{equation}\label{DellDef}
D_\ell\ := \left[\;
\parbox{270pt}{lowest common denominator of the $\Psi_J$ for all $J$ with $\ell = \bigl|\{i\in J\mid i\le s\}\bigr|$ satisfying Condition \ref{1JConds} (a).}\;\right]
\end{equation}

\begin{Lem}\label{1LCD0}
We have $D_0=1-Y^{pq}$.
\end{Lem}

\begin{Proof}
Consider any subset $J\subset\{s+1,\ldots,r\}$. Then $J$ trivially satisfies Condition \ref{1JConds} (a). 
Also, for all $m\ge0$ we have $1\not\in\pi^m(J)$ and $\ell_m=0$.
Thus the first sum in (\ref{1PsiJDef}) is zero and second is the sum of $Y^{m+s+1}$ for all $m\ge0$ with $s+1\in\pi^m(J)$. Since the condition $s+1\in\pi^m(J)$ depends only on $m\bmod(r-s)$, we deduce that 
$$(1-Y^{r-s})\cdot \Psi_J \ = \hskip-5pt \sum_{0\le m<r-s \atop s+1\in\pi^m(J)} \kern-10pt Y^{m+s+1}.$$
Moreover, for any $0\le m<r-s$ we have $s+1\in\pi^m(J)$ if and only if $m+s+1=\pi^{-m}(s+1)$ lies in~$J$. Therefore
$$(1-Y^{r-s})\cdot \Psi_J \ = \sum_{s<j\in J} Y^j.\hskip26pt$$
In particular the denominator of $\Psi_J$ divides $(1-Y^{r-s})$. Conversely, in the case $J=\{r\}$ we obtain $(1-Y^{r-s})\cdot \Psi_{\{r\}} = Y^r$. Thus the lowest common denominator is $1-Y^{r-s}=1-Y^{pq}$, as desired.
\end{Proof}

\begin{Lem}\label{1LCDs}
We have $D_s=1-XY$.
\end{Lem}

\begin{Proof}
{}From Lemma \ref{1JCondsaEq} we see that the only subset $J\subset \{1,\ldots,r\}$ containing $\{1,\ldots,s\}$ which satisfies Condition \ref{1JConds} (a) is $\{1,\ldots,r\}$ itself. For it we have $1\in\pi^m(J)$ and $\ell_m=m$ for all $m\ge0$. By (\ref{1PsiJDef}) this implies that $\Psi_J=\sum_{m\ge0}X^mY^{m+1} = Y/(1-XY)$, whose denominator is $1-XY$, as desired.
\end{Proof}

\begin{Lem}\label{1LCDLemUpperAll}
For any $0<\ell<s$, the lowest common denominator of the $\Psi_J$ for all subsets $J\subset\{1,\ldots,r\}$ satisfying $\ell = \bigl|\{i\in J\mid i\le s\}\bigr|$ and  $\{s+1,\ldots,r\} \subset J$ is $1-X^{\ell}Y^{s}$.
\end{Lem}

\begin{Proof}
For any such $J$ we have $s+1\in\pi^m(J)$ for all $m\ge0$. Thus in both sums in (\ref{1PsiJDef}), the terms with $m+s$ in place of $m$ are obtained on multiplying the terms for $m$ by $X^\ell Y^s$. Therefore 
$$(1-X^{\ell}Y^{s})\Psi_J \ =  
\sum_{0\le m<s \atop 1\in\pi^m(J)} \!\! X^{\ell_m}Y^{m+1} \; + 
\!\!\sum_{0\le m<s \atop 1\not\in\pi^m(J)} \!\! X^{\ell_m}Y^{m+s+1}.$$
Thus the lowest common denominator divides $(1-X^{\ell}Y^{s})$. 

For the converse note first that for all ${0\le m<s}$ we have $1\not\in\pi^m(J)$ if and only if $m+1=\pi^{-m}(1)\in J$.  Thus for all $0\le m<s$ we have $\ell_m = |\{i\in J\mid i\le m\}$. 
Now consider the subsets 
$$\begin{array}{ll}
J  \!&:=\ \{2,3,\ldots,\ell+1,s+1,\ldots,r\}, \\[3pt]
J' \!&:=\ \{1,3,\ldots,\ell+1,s+1,\ldots,r\},
\end{array}$$
both of which satisfy the given conditions. The preceding remarks show that in the range $0\le m<s$, the values of $\ell_m$ associated to $J$ and $J'$ differ only for $m=1$. Thus the summands for all $m>1$ in both sums above are the same for $J$ and~$J'$, and the others yield
$$(1-X^{\ell}Y^{s})(\Psi_J-\Psi_{J'}) \ =\ 
(Y^2+Y^{s+1}) - (Y+XY^{s+2}).$$
As the right hand side is coprime to $(1-X^{\ell}Y^{s})$, the lemma follows.
\end{Proof}

\begin{Lem}\label{1LCDBig}
For any $s-\frac{s}{p}<\ell<s$ we have $D_\ell=1-X^\ell Y^s$.
\end{Lem}

\begin{Proof}
Since $p$ divides~$s$, any residue class modulo $p$ contains precisely $\frac{s}{p}$ elements from $\{1,\ldots,s\}$. Thus if $\{i\in J\mid i\le s\}$ has cardinality $\ell >s-\frac{s}{p}$, it must meet every residue class modulo~$p$. If in addition $J$ satisfies Condition \ref{1JConds} (a), then Lemma \ref{1JCondsaEq} implies that ${\{s+1,\ldots,r\}} \subset \nobreak J$. Conversely, any subset $J\subset\{1,\ldots,r\}$ with $\ell = \bigl|\{i\in J\mid i\le s\}\bigr|$ and $\{s+1,\ldots,r\} \subset J$ trivially satisfies Condition \ref{1JConds} (a). The lemma therefore reduces to Lemma \ref{1LCDLemUpperAll}.
\end{Proof}

\begin{Lem}\label{1LCDRedLem}
Consider any subset $J\subset\{1,\ldots,r\}$ with $\ell = \bigl|\{i\in J\mid i\le s\}\bigr|$ that satisfies Condition \ref{1JConds} (a).
Then $\bar J := J\cup\{s+1,\ldots,r\}$ has the same properties and
$$\Psi_{\bar J} - \Psi_J \ =
 \kern-14pt \sum_{m\ge0  \atop s+1\in\pi^m(\bar J\setminus J)}
\kern-15pt X^{\ell_m}Y^{m+s+1}.$$
\end{Lem}

\begin{Proof}
The statement about $\bar J$ follows from the form of Condition \ref{1JConds} (a). Next, the condition $1\in\pi^m(J)$ and the exponent $\ell_m$ in (\ref{1PsiJDef}) is the same for $J$ as for~$\bar J$. Thus the difference $\Psi_{\bar J} - \Psi_J$ comes only from the terms of the second sum in (\ref{1PsiJDef}) with $1\not\in\pi^m(J)$ and $s+1\in\pi^m(\bar J\setminus J)$. But since $J$ satisfies Condition \ref{1JConds} (a), any $m$ with $s+1\in\pi^m(\bar J\setminus J)$ already satisfies $1\not\in\pi^m(J)$. Thus the indicated formula follows.
\end{Proof}

\begin{Lem}\label{1LCDSmall}
For any $0<\ell<s-\frac{s}{p}$ we have $D_\ell=1-X^{\ell q}Y^{sq}$.
\end{Lem}

\begin{Proof}
Lemma \ref{1PsiRat} already shows that the lowest common denominator $D_\ell$ divides ${1-X^{\ell q}Y^{sq}}$. For the converse we apply Lemma \ref{1LCDRedLem} to the case that ${J\cap\{s+1,\ldots,r\}} \allowbreak = {\{s+2,\ldots,r\}}$. Then the condition $s+1\in\pi^m(\bar J\setminus J) = \{\pi^m(s+1)\}$ is equivalent to $r-s=pq\mid m$. Moreover, as in the proof of Lemma \ref{1PsiRat}, the terms of the sum with $m+sq$ in place of $m$ are obtained on multiplying the terms for $m$ by $X^{\ell q}Y^{sq}$. Therefore 
$$(1-X^{\ell q}Y^{sq})\cdot(\Psi_{\bar J} - \Psi_J) \ =
 \sum_{0\le m<sq \atop pq|m} \kern-5pt X^{\ell_m}Y^{m+s+1}.$$
Thus it suffices to show that some linear combination of this for all possible $J$ is coprime to $1-X^{\ell q}Y^{sq}$. In fact, the difference for two suitable choices of $J$ will do.

Recall that $p=\gcd(s,r-s)$, so that $s-\frac{s}{p}$ is the cardinality of $\{1\le i\le s: p\ndiv i\}$. Also note that the assumption $0<\ell<s-\frac{s}{p}$ implies that $p=\gcd(s,r-s)>1$. Thus we can choose a subset $A$ of $\{1\le i\le s: p\ndiv i\}$ of cardinality $\ell$ such that $1\not\in A$ and $p+1\in A$. Then $A' := \{1\}\cup A\setminus\{p+1\}$ is another subset of $\{1\le i\le s: p\ndiv i\}$ of cardinality~$\ell$. With these choices we set $J := A\cup\{s+2,\ldots,r\}$ and $J' := A'\cup\{s+2,\ldots,r\}$. Then in view of Lemma \ref{1JCondsaEq}, both $J$ and $J'$ satisfy the stated conditions.

Next recall that the exponent $\ell_m$ for $J$ was defined as $\ell_m := |\{0\le i<m\mid 1\in\pi^i(J)\}|$. Define accordingly $\ell'_m := |\{0\le i<m\mid 1\in\pi^i(J')\}|$. 
Since $\pi$ induces a permutation of order $s$ on $\{1,\ldots,s\}$, these numbers satisfy $\ell_{m+s}=\ell_m+\ell$ and $\ell'_{m+s}=\ell'_m+\ell$. Thus the difference $\ell'_m-\ell_m$ depends only on $m\bmod(s)$.

Recall also that for all ${0\le m<s}$ we have 
$\ell_m = |\{i\in J\mid i\le m\}$, and similarly $\ell'_m = 
|\{i\in J'\mid i\le m\}|$. The construction of $J$ and $J'$ thus implies that these values are equal unless $1\le m\le p$, in which case $\ell'_m = \ell_m+1$.

Returning now to the above sum, consider any integer $0\le m<sq$ with $pq|m$. Write $m=np+ks$ with $0\le np<s$. Then the preceding remarks imply that $\ell'_m-\ell_m = \ell'_{np}-\ell_{np}$ is $0$ unless $n=1$, in which case it is~$1$. The latter case occurs if and only if $m\equiv p\bmod(s)$. 
But by the definition of~$p$, the integers $\frac{s}{p}$ and $\frac{r-s}{p}=q$ are relatively prime. Thus there is precisely one integer $0\le m<sq$ with $pq|m$ and $m\equiv p\bmod(s)$. For this $m$ we then have
$$(1-X^{\ell q}Y^{sq})\cdot(\Psi_{J\cup\{s+1\}} - \Psi_J
-\Psi_{J'\cup\{s+1\}} + \Psi_{J'}) \ =
\ X^{\ell_m}(1-X)Y^{m+s+1}.$$
As the right hand side is coprime to ${(1-X^{\ell q}Y^{sq})}$, the lemma follows.
\end{Proof}

\begin{Lem}\label{1LCDLemMid}
For $\ell=s-\frac{s}{p}>0$ we have 
$$D_\ell \ =\ \frac{(1-X^{\ell}Y^{s})(1-X^{pq-q}Y^{pq})}{(1-X^{p-1}Y^{p})}.$$
\end{Lem}

\begin{Proof}
Consider any subset $J\subset\{1,\ldots,r\}$ with $\ell = \bigl|\{i\in J\mid i\le s\}\bigr|$ that satisfies Condition \ref{1JConds} (a). Set $\bar J := J\cup\{s+1,\ldots,r\}$, which again has the stated properties. Then by Lemma \ref{1LCDLemUpperAll} the lowest common denominator of $\Psi_{\bar J}$ for all possible $J$ is $1-X^{\ell}Y^{s}$. 

Suppose in addition that $J\not=\bar J$, and choose an element $j_0\in\bar J\setminus J$. Then Lemma \ref{1JCondsaEq} implies that 
$$\{i\in J\mid i\le s\} \ \subset\ \{1\le i\le s\mid i\not\equiv j_0\bmod(p)\}.$$
Since $p$ divides~$s$, the set on the right hand side has cardinality $s-\frac{s}{p}=\ell$, same as the set on the left hand side. Thus the inclusion is an equality, and so
\UseTheoremCounterForNextEquation
\begin{equation}\label{1LCDLemEq1}
\{1\le i\le r\mid i\not\equiv j_0\bmod(p)\}
\ \subset\ J \ \subset\ \{1\le i\le r\mid i\not\equiv j_0\bmod(p) \;\lor\; i>s\}.
\end{equation}
In particular this implies that $\ell_{m+p}=\ell_m+p-1$ for all~$m\ge0$. Moreover, the condition $s+1\in\pi^m(\bar J\setminus J)$ depends only on $m\bmod(r-s)$ with $r-s=pq$. It follows that in the sum in Lemma \ref{1LCDRedLem}, the terms with $m+pq$ in place of $m$ are obtained on multiplying the terms for $m$ by $X^{pq-p}Y^{pq}$. Therefore 
\UseTheoremCounterForNextEquation
\begin{equation}\label{1LCDLemEq2}
(1-X^{pq-q}Y^{pq})\cdot(\Psi_{\bar J} - \Psi_J) \ =
 \kern-14pt \sum_{0\le m<pq \atop s+1\in\pi^m(\bar J\setminus J)}
\kern-15pt X^{\ell_m}Y^{m+s+1}.
\end{equation}
In particular, the denominator of $\Psi_{\bar J} - \Psi_J$ divides $1-X^{pq-q}Y^{pq}$.

Conversely, the subset $J := \{1\le i\le r\mid i\not\equiv 1\bmod(p)\;\lor\; i>s+1\}$ satisfies (\ref{1LCDLemEq1}) with $j_0:= s+1$. For it we have $\bar J\setminus J = \{s+1\}$, and so the conditions $0\le m<pq$ and $s+1\in\pi^m(\bar J\setminus J)$ in (\ref{1LCDLemEq2}) are satisfied only for $m=0$. The right hand side of (\ref{1LCDLemEq2}) is therefore equal to $Y^{s+1}$ in this case. Together it follows that the lowest common denominator of $\Psi_{\bar J}-\Psi_J$ for all possible $J$ is $1-X^{pq-q}Y^{pq}$. 

Combining everything, we deduce that $D_\ell$ is the least common multiple of $1-X^{\ell}Y^{s}$ and $1-X^{pq-q}Y^{pq}$. But since $\ell=(p-1)\frac{s}{p}$ and, by the definition of~$p$, the integers $\frac{s}{p}$ and $q$ are relatively prime, the greatest common divisor of $1-X^{\ell}Y^{s}$ and $1-X^{pq-q}Y^{pq}$ is $1-X^{p-1}Y^{p}$. Thus the least common multiple has the indicated form.
\end{Proof}

\medskip
Combining the preceding lemmas, we obtain:

\begin{Prop}\label{1LCDAll}
For any $0\le\ell\le s$ we have
$$D_{\ell} \ =\ 
\left\{\begin{array}{ll}
(1-Y^{pq}) & \hbox{if $\ell=0$,} \\[5pt]
(1-X^{\ell q}Y^{sq}) & \hbox{if $0<\ell<s-\frac{s}{p}$,} \\[5pt]
\displaystyle\frac{(1-X^{\ell}Y^{s})(1-X^{pq-q}Y^{pq})}{(1-X^{p-1}Y^{p})}
& \hbox{if $\ell=s-\frac{s}{p}>0$,} \\[13pt]
(1-X^{\ell}Y^{s}) & \hbox{if $s-\frac{s}{p}<\ell<s$,} \\[5pt]
(1-XY) & \hbox{if $\ell=s$.} \\
\end{array}\right.$$
\end{Prop}


\subsection{Denominators of orbit length generating functions}
\label{sec1Denom}

\begin{Lem}\label{1BiOLGFLem}
For all $1\le i\le r$ we have 
$$\Psi_{b_i} \ =\ 
\left\{\begin{array}{cl}
Y^i & \hbox{if $i\le s$,}\\[3pt]
\displaystyle \frac{Y^i}{1-Y^{r-s}} & \hbox{if $i>s$.}
\end{array}\right.$$
\end{Lem}

\begin{Proof}
By the recursion relations (\ref{1bRecRels}) and Proposition \ref{OLGFRecRels22}, and the fact that $\Psi_1=0$ by (\ref{F12}), we have
$$\Psi_{b_i} \ =\ 
\left\{\begin{array}{ll}
Y & \hbox{if $i=1$,}\\[3pt]
Y\Psi_{b_r}+Y\Psi_{b_s} & \hbox{if $i=s+1$,}\\[3pt]
Y\Psi_{b_{i-1}} & \hbox{otherwise.}
\end{array}\right.$$
By induction on $i$ this implies that $\Psi_{b_i}=Y^i$ for all $1\le i\le s$. Induction also shows that $\Psi_{b_i}=Y^{i-s-1}\Psi_{b_{s+1}}$ for all $s<i\le r$. Therefore
$\Psi_{b_{s+1}} = Y\Psi_{b_r}+Y\Psi_{b_s} = Y^{r-s}\Psi_{b_{s+1}}+Y^{s+1}$ and hence $\Psi_{b_{s+1}} = Y^{s+1}/(1-Y^{r-s})$.
This in turn implies that $\Psi_{b_i} = Y^i/(1-Y^{r-s})$ for all $s<i\le r$, and we are done.
\end{Proof}

\begin{Prop}\label{1IrredDenom}
For any \OLGF-irreducible element $w\in\GammaUx$ we have $\Phi_w = \Phi_{J_w}$.
\end{Prop}

\begin{Proof}
We apply the constructions of Subsection \ref{sec1Prop} to $w_0 := w$. Combining Lemma \ref{1LemInduc} for all $n\ge0$ with Proposition \ref{OLGFRecRels22} and the conjugation invariance of orbit length generating functions yields
$$\Psi_{w_n} \ =\ 
\left\{\begin{array}{ll}
Y+XY\Psi_{w_{n+1}} & \hbox{if $1\in J_{w_n}$,} \\[3pt]
Y\Psi_{b_s}+Y\Psi_{w_{n+1}} & \hbox{if $1\not\in J_{w_n}$ and $s+1\in J_{w_n}$,} \\[3pt]
Y\Psi_1+Y\Psi_{w_{n+1}} & \hbox{if $1\not\in J_{w_n}$ and $s+1\not\in J_{w_n}$.}
\end{array}\right.$$
Here $\Psi_1=0$ by (\ref{F12}), and $\Psi_{b_s}=Y^s$ by Lemma \ref{1BiOLGFLem}. Moreover, by Lemma \ref{1JIter} we have $J_{w_n}=\pi^n(J_{w_0})$ for all $n\ge0$. As in (\ref{1PsiJDef}) let $\ell_m$ denote the number of integers $0\le \ell<m$ such that $1\in\pi^\ell(J_{w_0})$. 
Then by induction on $n$ it follows that 
$$\Psi_{w_0} \ =\ 
\! \sum_{0\le m<n \atop 1\in\pi^m(J_{w_0})} \kern-10pt X^{\ell_m}Y^{m+1} \;+ \kern-13pt \sum_{0\le m<n \atop 1\not\in\pi^m(J_{w_0})\ni s+1} 
\kern-16pt X^{\ell_m}Y^{m+s+1}
\; + \; X^{\ell_n}Y^n\Psi_{w_n}$$
In the limit over $n$ this implies that $\Psi_{w_0} = \Psi_{J_{w_0}}$. Using (\ref{1PhiJDef}) and Proposition \ref{OLGFRecRels21} we deduce that $\Phi_{w_0} = \Phi_{J_{w_0}}$, as desired.
\end{Proof}


\medskip
For any $0\le\ell\le s$ we now define
\UseTheoremCounterForNextEquation
\begin{equation}\label{DellUxDef}
D_{\Ux,\ell}\ := \left[\;
\parbox{270pt}{lowest common denominator of the $\Psi_J$ for all $J$ with $\ell = \bigl|\{i\in J\mid i\le s\}\bigr|$ satisfying Conditions \ref{1JConds}.}\;\right]
\end{equation}
By construction this is a divisor of the polynomial $D_\ell$ from (\ref{DellDef}) and Proposition \ref{1LCDAll}.

\begin{Thm}\label{1MainThm}
The power series $\Phi_w \in 1+Y\BZ[[X,Y]]$ for all $w\in\GammaUx$ are rational functions with the lowest common denominator 
$$D_\Ux\ :=\ (1-2Y)\;\cdot\! \prod_{0\le\ell\le s}\!D_{\Ux,\ell}
\ \in\ 1+Y\BZ[X,Y].$$
\end{Thm}

\begin{Proof}
By Proposition \ref{1OLGFFin} and Theorem \ref{OLGFLinComb} the $\Phi_w$ for all $w\in\GammaUx$ are $\BZ[X,Y]$-linear combinations of the $\Phi_w$ for all \OLGF-irreducible elements. By 
Propositions \ref{1JCondLem} and \ref{1IrredExists} and \ref{1IrredDenom} the latter are precisely the $\Phi_J$ for all subsets $J\subset\{1,\ldots,r\}$ which satisfy Conditions \ref{1JConds}. They are therefore rational functions, and in view of (\ref{1PhiJDef}) and (\ref{DellUxDef}) their lowest common denominator is the least common multiple of the polynomials $(1-2Y)D_{\Ux,\ell}$ for all~$\ell$. 

Lemma \ref{1PsiRat} shows that each $D_{\Ux,\ell}$ divides $1-X^{\ell q}Y^{sq}$. Thus $D_{\Ux,\ell}$ can be chosen congruent to $1\bmod Y$, and then $D_\Ux$ has the same property. 
Moreover, the polynomials $1-2Y$ and $1-X^{\ell q}Y^{sq}$ for all $0\le\ell\le s$ are pairwise coprime, for instance because, viewed as polynomials in~$Y$, their zeros in an algebraic closure of $\BQ(X)$ are mutually distinct. Thus the least common multiple of all $(1-2Y)D_{\Ux,\ell}$ is~$D_\Ux$, and we are done.
\end{Proof}

\begin{Prop}\label{1Lem9}
\begin{enumerate}
\item[(a)] For all $\Ux$ we have $D_{\Ux,0}=D_0=1-Y^{pq}$ and $D_{\Ux,s}=D_s=1-XY$.
\item[(b)] For $\Ux=(0,\ldots,0)$ or $(1,\ldots,1)$ we have $D_{\Ux,\ell}=D_\ell$ for all~$\ell$.
\end{enumerate}
\end{Prop}

\begin{Proof}
The first statement in (a) is a consequence of Lemmas \ref{1LCD0} and \ref{1BiOLGFLem}. The second statement in (a) results from Lemma \ref{1LCDs} and the fact that $J:= \{1,\ldots,r\}$ always satisfies Conditions \ref{1JConds}. Assertion (b) is a consequence of the fact that for $\Ux=(0,\ldots,0)$ or $(1,\ldots,1)$ Condition \ref{1JConds} (b) is satisfied for all~$J$.
\end{Proof}


\medskip
In principle, the determination of the lowest common denominator $D_\Ux$ in Theorem \ref{1MainThm} is a finite combinatorial problem concerning the data $r$, $s$, and~$\Ux$. The author does not (yet) know a simple direct description in general. However, we determined $D_\Ux$ in small cases using the computer algebra system Maple: see \cite{Pink2014Maple}. The outcome was that whenever $r\le 8$ and $x_2,\ldots,x_r$ are not all equal, then $\prod_{0<\ell<s}D_{\Ux,\ell}=1$ except in the following cases:
$$\begin{array}{|c|l|l|}
\hline
{\Large\strut} (r,s) & \, \prod_{0<\ell<s}D_{\Ux,\ell} & \hbox{Conditions on $\Ux=(x_2,\ldots,x_r)$} \\[5pt]
\hline\hline
{\large\strut} (4,2) & (1{-}XY^2) & x_2=x_4 \\
\hline
{\large\strut} (5,4) & (1{-}XY^2) & x_2=x_4=x_5\not=x_3 \\
\hline
{\large\strut} (6,2) & (1{-}XY^2) & x_2=x_4=x_6 \hbox{\ but not all $x_i$ equal} \\
\hline
{\large\strut} (6,3) & (1{-}XY^3)(1{-}X^2Y^3) & x_2=x_3=x_5=x_6 \\[-2pt]
{\large\strut} (6,3) & (1{-}XY^3) & x_2=x_5\not=x_3=x_6 \\
\hline
{\large\strut} (6,4) & (1{-}XY^4)(1{-}XY^2) & x_2=x_4=x_6\not=x_3=x_5 \\[-2pt]
{\large\strut} (6,4) & (1{-}XY^2) & x_2=x_4=x_6 \longland x_3\not=x_5 \\
\hline
{\large\strut} (7,4) & (1{-}XY^2) & x_2=x_4=x_5=x_6=x_7\not=x_3 \\
\hline
{\large\strut} (7,6) & (1{-}XY^3)(1{-}X^2Y^3) & x_2=x_3=x_5=x_6=x_7\not=x_4 \\[-2pt]
{\large\strut} (7,6) & (1{-}XY^2) & x_2=x_4=x_6=x_7 \hbox{\ but not all $x_i$ equal} \\
\hline
{\large\strut} (8,2) & (1{-}XY^2) & x_2=x_4=x_6=x_8 \hbox{\ but not all $x_i$ equal} \\
\hline
{\large\strut} (8,4) & (1{-}XY^4)(1{-}X^2Y^4)(1{-}X^3Y^4) & x_2=x_3=x_4=x_6=x_7=x_8 \\[-2pt]
{\large\strut} (8,4) & (1{-}XY^4)(1{-}XY^2) & x_2=x_4=x_6=x_8\not=x_3=x_7 \\[-2pt]
{\large\strut} (8,4) & (1{-}X^2Y^4) & x_2=x_4=x_5=x_6=x_7=x_8\not=x_3 \\[-2pt]
{\large\strut} (8,4) & (1{-}XY^4) & x_2=x_6\not=x_4=x_8 \longland x_3=x_7 \\[-2pt]
{\large\strut} (8,4) & (1{-}XY^2) & x_2=x_4=x_6=x_8 \longland x_3\not=x_7 \land (x_5,x_7)\not=(0,0)\! \\
\hline
{\large\strut} (8,6) & (1{-}XY^6)(1{-}X^2Y^6)(1{-}XY^2) & x_2=x_4=x_6=x_8\not=x_3=x_5=x_7 \\[-2pt]
{\large\strut} (8,6) & (1{-}XY^3)(1{-}X^2Y^3) & x_2=x_3=x_5=x_6=x_7=x_8\not=x_4 \\[-2pt]
{\large\strut} (8,6) & (1{-}XY^2) & x_2=x_4=x_6=x_8 \hbox{\ but $x_3$, $x_5$, $x_7$ not all equal} \\
\hline
\end{array}$$

%
%
%
%
%
%
%
%

 \newpage

%



\end{document}